\documentclass[english]{elsarticle}
\usepackage[english]{babel}
\usepackage[utf8]{inputenc}
\usepackage{geometry}
\usepackage{yfonts}
\usepackage{xfrac}
\usepackage[c]{esvect}
\usepackage{setspace}
\usepackage{graphicx}
\usepackage{indentfirst} 
\usepackage{amsmath,amsfonts,amsthm,amscd,upref,amstext}	
\usepackage[titletoc]{appendix}
\usepackage[utf8]{inputenc} 
\usepackage[T1]{fontenc}    
\usepackage{hyperref}       
\usepackage{url}            
\usepackage{booktabs}       
\usepackage{amsfonts}       
\usepackage{amsmath}
\usepackage{nicefrac}       
\usepackage{microtype}      
\usepackage{amsthm}
\usepackage{graphicx}       
\usepackage{caption}
\usepackage{comment}
\usepackage{tikz}
\usepackage{tikz-cd}
\usepackage{thmtools, thm-restate}
\usepackage{cleveref}

\newtheorem{theorem}{Theorem}[section]
\newtheorem{lemma}{Lemma}[section]
\newtheorem{definition}{Definition}[section]
\numberwithin{equation}{section}

\title{Affine Super Yangian and Weyl groupoid}

\author{Vladimir Stukopin}
\address{Moscow Institute of Physics and Technology, Center of Fundamental Mathematics,\\  South Mathematical Institute}

\author{Vasiliy Volkov}
\address{Moscow Institute of Physics and Technology, Center of Fundamental Mathematics}

\begin{document}

\begin{abstract}
We define affine Super Yangian $Y_{\hbar}(\hat{sl}(m|n), \Pi) $ for affine special linear superalgebra $\hat{sl}(m|n)$ and arbitrary system of simple roots $\Pi$ in terms of minimalistic system of generators. We also consider Drinfeld presentation for affine super Yangian in the case of arbitrary simple root system $\Pi$ and prove that these two presentations (Drinfeld and minimalistic) of $Y_{\hbar}(\hat{sl}(m|n), \Pi)$ are isomorphic as associative superalgebras.  We also construct isomorphism of affine super Yangians $Y_{\hbar}(\hat{sl}(m|n), \Pi)$ and $Y_{\hbar}(\hat{sl}(m|n), \Pi')$ for different simple root systems $\Pi$ and $\Pi'$. After them we also define Weyl groupoid as a set of morphisms in category with objects, which are super Yanginas $Y_{\hbar}(\hat{sl}(m|n), \Pi)$, where $\Pi$ is simple root system. We describe Weyl groupoid in terms of generators and describe action of these generators on super Yangians. We describe isomorphisms between $Y_{\hbar}(\hat{sl}(m|n), \Pi)$ and $Y_{\hbar}(\hat{sl}(m|n), \Pi')$ as elements of Weyl groupoid.
\end{abstract}



\maketitle


\section{Introduction}
We define affine Super Yangian $Y_{\hbar}(\hat{sl}(m|n), \Pi)$ for affine special linear superalgebra $\hat{sl}(m|n)$ for arbitrary system of simple roots $\Pi$ in terms of minimalistic system of generators. We also consider Drinfeld presentation of $Y_{\hbar}(\hat{sl}(m|n), \Pi)$ and constructed isomorphism of these two presentations of  affine Yangian for  arbitrary simple root system $\Pi$. After them  we define affine Weyl groupoid, and define action of affine Weyl groupoid on different Yangian realizations. We show that action of Weyl groupoid sets the isomorphisms of Yangian realizations, defined by different systems of simple root systems.

The construction of the Yangian in so called Faddeev-Reshetikhin-Takhtadjan (FRT) presentation appeared in connection with the application of the algebraic Bethe ansatz to the study of quantum integrable models with a rational $R$-matrix before the appearance of the term itself. The concept of the Yangian was given by Drinfeld as one of the most important for applications example of quantum groups (\cite{Dr}). Drinfeld defined the Yangian of a finite dimensional simple Lie algebra $\mathfrak{g}$ in order to obtain a solution of the Yang-Baxter equation. The Yangian is a quantum group which is the deformation of the current algebra $\mathfrak{g}[z]$. He defined it by three different presentations. One of those presentations is called the Drinfeld presentation (\cite{Drinfeld}) (or new Drinfeld presentation) whose generators are $\{h_{i,r}, x^{\pm}_{i,r}| r\in \mathbb{Z}_{\geq 0}\}$, where $\{h_i, x^{\pm}_i\}$ are Chevalley generators of $\mathfrak{g}$. Yangian is a Hopf algebra and has comultiplication.

The definition of Yangian as an associative algebra naturally extends to the case that $\mathfrak{g}$ is a symmetrizable affine Kac-Moody Lie algebra in the Drinfeld presentation \cite{Drinfeld}. Such type of a Yangian is called affine Yangian and it was first determined by apparently S. Boyarchenko and S. Levendorskii \cite{Levendorski} (see also work \cite{Guaywork}).

Defining its quasi-Hopf algebra structure is more involved, but this problem has been settled for Yangian of  affine Kac-Moody Lie algebras in works by N.Guay with collaborators. It is known that the Yangians are closely related to $W$-algebras. It was shown  that there exist surjective homomorphisms from Yangians of type $A$ to rectangular finite $W$-algebras of type $A$. Affine Yangian related to infinite $W$-algebras which play important role in mathematical physics.
In the case of the Lie superalgebra $sl(m|n)$, also it is known definition of Yangian as in the Drinfeld presentation as so called FRT presentation. The relationship between Yangians and $W$-algebras were also studied in the case of finite Lie superalgbras by many authors: E. Ragoucy, P. Sorba (\cite{RS}),    C. Peng (\cite{Peng}), V. Serganova and E. Poletaeva (\cite{PS}). In the recent paper R. Gaberdiel, W. Li, C. Peng and H. Zhang \cite{SuSy} defined the Yangian $Y(\hat{gl}(1|1))$ for the affine Lie superalgebra $\hat{gl}(1|1)$. Ueda defined affine super Yangian $Y(\hat{sl}(m|n))$ for distinguished simple root system.

It should be noted that the basic Lie superalgebra, in contrast to the simple Lie algebra, can be given by different Dynkin diagrams, what is explained by the fact that it has different, nonequivalent systems of simple roots (or, which is the same, has nonconjugate Borel subalgebras). We define super Yangian $Y_{\hbar}(\hat{sl}(m|n))$ for arbitrary simple root system $\Pi$. We also consider the super Yangian $Y_{\hbar}(\tilde{sl}(m,n))$, where $\hat{sl}(m|n)$ is a presentation of affine Kac-Moody algebra as central extension of double loop superalgebra. We prove that for any two different simple root systems $\Pi_1$ and $\Pi_2$ the corresponding affine super Yangians $Y_{\hbar}(\tilde{sl}(\Pi_1))$ and $Y_{\hbar}(\tilde{sl}(\Pi_2))$ are isomorphic. The same result can be proved also for super Yangians $Y_{\epsilon_1, \epsilon_2}(\tilde{sl}(\Pi_1))$ and $Y_{\epsilon_1, \epsilon_2}(\tilde{sl}(\Pi_2))$  (which are introduced in \cite{Guaywork1} for affine Yangian of $\hat{sl}(n)$ and close related to quantum toroidal superalgebras).  Our proof is based on construction of affine Weyl groupoid which generated by (super)reflections of weight lattice relatively simple roots. The (super)reflections induce above mentioned isomorphisms.  We construct two presentations of affine  super Yangian, namely, so called minimalistic presentation  and Drinfeld presentation. We note that in the case Yangian such presentation was introduced by S. Levendorskii \cite{Lev} and for super Yangians minimalistic presentation was done in \cite{St}.  Our second result is that this two presentations are isomorphic as associative superalgebras. We will not consider questions related to comultiplication here. This is the topic of a separate work, which is planned as a continuation of this study.

\section{Preliminaries}
\label{sec:headings}
\subsection{Basic Definitions}

\paragraph{Classical Lie superalgebras}

\begin{definition}  {\it Lie Superalgebra} is a $\mathbb{Z}_2$-graded vector space $\mathfrak{g}=\mathfrak{g}_0\oplus \mathfrak{g}_1$ together with a bilinear map $[,]:\mathfrak{g}\times \mathfrak{g} \rightarrow \mathfrak{g}$ and parity function $p$ ($p(x) = i$ if $x \in g_i$) such that
\begin{enumerate}
    \item $[\mathfrak{g}_\alpha,\mathfrak{g}_\beta]\subseteq \mathfrak{g}_{\alpha+\beta} \quad  \alpha, \beta \in \mathbb{Z}_2$ \quad ($\mathbb{Z}_2$-grading),
    \item $[a,b]=-(-1)^{\bar{a}\bar{b}}[b,a]$ \quad (graded skew-symmetry),
    \item $(-1)^{\bar{a}\bar{c}}[a,[b,c]]+(-1)^{\bar{a}\bar{b}}[b,[c,a]]+(-1)^{\bar{b}\bar{c}}[c,[a,b]=0$ \quad (graded Jacobi identity),
\end{enumerate}
for all $a,b,c \in g$ (see \cite{Musson}).
\end{definition}
We write $\bar{a}=p(a)$ if $a \in \mathfrak{g}_{\bar{a}}$.
Important example of Lie superalgebras are classical Lie superalgebras. Let $g$ be finite dimensional Lie superalgebra such that $\mathfrak{g}_0$ is reductive and $\mathfrak{g}_1$ is a semisimple $\mathfrak{g}_0$-module. Let $\mathfrak{h}_0$ be a Cartan subalgebra of $\mathfrak{g}_0$. For $\alpha \in \mathfrak{h}_{0}^{\ast}$ set
\[
\mathfrak{g}^{\alpha}=
\begin{cases}
x \in \mathfrak{g}|[h,x]=\alpha(h)x\quad \forall h \in \mathfrak{h}_0\}
\end{cases}
\]
and let
\[
\Delta=
\begin{cases}
\alpha \in \mathfrak{h}_0^{\ast}|\alpha \neq 0,\mathfrak{g}^{\alpha} \neq 0\}
\end{cases}
\]
be the set of roots $\mathfrak{g}$.Since the action of $\mathfrak{h}_0$ on any finite dimensional simple $\mathfrak{g}$-module is diagonalizable. Thus there {\it root space decomposition}
\begin{equation}
\label{eq:5}
\mathfrak{g}=\mathfrak{h}\oplus_{\alpha \in \Delta} \mathfrak{g}^{\alpha},
\end{equation}
where $\mathfrak{h} = \mathfrak{g}^0$ is the centralizer of $\mathfrak{h}_0$ in $\mathfrak{g}$. Then we have following lemma
\begin{lemma}
 If $\mathfrak{g}$ is a classical simple Lie superalgebra and $\alpha,\beta,\alpha+\beta$ are roots of $\mathfrak{g}$, then $[\mathfrak{g}^{\alpha},\mathfrak{g}^{\beta}]=\mathfrak{g}^{\alpha+\beta}$. \\
\end{lemma}

Next we recall some definitions from Lie superalgebra theory related to the root systems.

\paragraph{Definition of root system}

\begin{definition}
A {\it Euclidean space} is a finite dimensional vector space $E$ equipped with positive definite inner product, denoted $(,)$.
\end{definition}
For nonzero $\alpha \in E$ we set $\alpha^{\vee}=2\alpha/(\alpha,\alpha)$. The  {\it hyperplane orthogonal} to $\alpha$ is
\[
\mathbb{H}_{\alpha}=\{\lambda \in E|(\lambda,\alpha)=0\}
\]
and we define the {\it reflection} $s_{\alpha}$ in $\mathbb{H}_{\alpha}$ by
\[
s_{\alpha}(\lambda)=\lambda-(\lambda,\alpha^{\vee})\alpha.
\]
Further we confine ourselves to the consideration of basic Lie superalgebras and the definitions will apply only to them.

\begin{definition}
An {\it irreducible root system} is a pair $(E,\Delta)$ where $E$ is a Euclidean space and $\Delta$ is a subset of $E$ such that the following axioms are satisfied.
\begin{enumerate}
    \item $\Delta$ is finite, spans $E$, and does not contain 0.
    \item If $\alpha \in \Delta$, then the only multiples of $\alpha$ contained in $\Delta$ are $\pm \alpha$.
    \item If $\alpha \in \Delta$, the reflection $s_{\alpha}$ leaves $\Delta$ invariant.
    \item If $\alpha, \beta \in \Delta$, then $(\alpha,\beta^{\vee})\in \mathbb{Z}$.
    \item $\Pi$ cannot be partitioned into two nonempty subsets $\Delta^{\prime}$ and $\Delta^{\prime \prime}$ such that $(\alpha,\beta)=0$ for all $\alpha \in \Delta^{\prime}$ and $\beta \in \Delta^{\prime \prime}$.
\end{enumerate}
\end{definition}

A subset $\Pi$ of $\Delta$ is called {\it basis of simple roots}, if $\Pi$ is a vector space basis of $E$ and each root $\beta$ of $\Phi$ can be rewritten in the form $\beta=\sum\limits_{\alpha \in \Pi}k_{\alpha}\alpha$ with integer coefficients $k_{\alpha}$, that are either all nonnegative or nonpositive. Then we suppose that $E=E_{1}\oplus E_{2}$ super vector space with parity function $p$, $p(v)=i$ if $v\in E_{i}$ the root $\beta=\sum\limits_{\alpha \in \Pi}k_{\alpha}\alpha$ is {\it positive} if all $k_{\alpha}$ are nonnegative.
The set of all positive roots $\beta$ we denote as $\Phi^{+}$. We denote set of even (odd) roots as $\Delta_{0}$($\Delta_{1}$). We will also formally consider the reflections related to odd roots (or superreflections).

\begin{definition}
Suppose that $(E,\Pi)$ is a root system. The {\it Weyl group} is the group generated by all reflections $s_{\alpha}$ with $\alpha \in \Pi^{+}$.\\
\end{definition}

\paragraph{Special linear superalgebra}

If $\mathfrak{g}$ is a Lie superalgebra of type $sl(m|n)=A(m-1,n-1)$, then $\mathfrak{g}$ has a root decomposition (\ref{eq:5}), where $\mathfrak{h}$ is set of diagonal matrices and $\Delta=\Delta_{0}\cup \Delta_{1}\subseteq \mathfrak{h}^{\star}$ set of roots. Let $\varepsilon_i,\delta_j$ be the linear functionals on $h$ whose values on the diagonal matrix
\[
a = diag(a_1,\ldots,a_{m+n+2})
\]
are given by
\begin{equation}
\label{rel:4}
\varepsilon_i(a)=a_i, \quad \delta_j(a)=a_{m+j}, \quad 1\leq i \leq m, \quad 1\leq j \leq n.
\end{equation}
Then
\begin{equation}
\Delta_0=\{\varepsilon_i-\varepsilon_j;\delta_i-\delta_j\}_{i\neq j}, \quad \Delta_1=\{\pm(\varepsilon_i-\delta_j)\}.
\end{equation}
We define Dynkin diagram as the graph with the set of nodes of two colours. Each node of white colour corresponds root from set $\Delta_0$, each node of grey colour corresponds root from set $\Delta_1$. Two nodes of the graph connected if they have nonzero dote product.
Note that we can introduce the function of order on system of simple roots $O(\Delta)\rightarrow \{1,\ldots, m+n+1\}$ by this definition, corresponding to Dynkin diagram.

{\it Even roots} is the elements of set $\Delta_0$, {\it odd roots} is the elements of set $\Delta_1$.\\

Note that inner (product) of two root from set $\Delta_1$ with itself equals $0$: $(\alpha_i,\alpha_i)=0$ for all roots from $\Delta_1$, nevertheless we can define reflection $s_{\alpha_i}$ for odd simple roots:
\begin{equation}
\label{rel:5}
    s_{\alpha_i}(\lambda)= \lambda + \alpha_i \quad if \quad \alpha+\lambda\quad  \text{is a root}, \quad
    s_{\alpha_i}(\alpha_i)=-\alpha_i,
    s_{\alpha_i}(\lambda)=\lambda, \quad otherwise.
\end{equation}
The set $s_{\alpha_i}(\Pi)=\{s_{\alpha_i}(\lambda)|\lambda \in \Pi\}$ there $\Pi$ is a basis of simple roots for $\mathfrak{g}$ is a basis of simple roots $\mathfrak{g}$.
Note that reflections have a groupoid structure on them, thus we can define {\it Weyl groupoid} as groupid generated by all reflections $s_{\alpha_i}$ for $\alpha_i \in \Delta$. The action of elements of Weyl groupoid generated by odd roots translates the system of simple roots $\Phi^+$ in isomorphic set of simple roots $\Delta^{+}_1$.

We need further the description of  the  basic Lie superalgbera $sl(m|n)$ as contragredient Lie superalgebra.

Then $sl(m|n)$ is isomorphic to the Lie super algebra  over $\mathbb{C}$ defined by generators $\{x_{i}^{\pm}, h_{i} | 1\leq i \leq m+n-1 \}$ and by the relations
\begin{equation}
\label{in:6}
[h_{i},h_{j}]=0,\quad    [h_i,x_{j}^{\pm}]=\pm a_{i,j}x_{j}^{\pm},\quad
[x_i^{+},x_j^{-}]=\delta_{i,j}h_i,\quad ad(x_{j}^{\pm})^{1+|a_{i,j}|}(x_{i}^{\pm}) = 0, \quad
\end{equation}
\begin{equation}
\label{in:7}
[x_{m}^{\pm},x_{m}^{\pm}]=0, \quad [[x_{m-1}^{\pm},x_{m}^{\pm}],[x_{m+1}^{\pm},x_{m}^{\pm}]]=0.
\end{equation}

\subsection{Affine Kac-Moody superalgebras}

We recall the definition of affinization $\hat{sl}(m|n)$ of $sl(m|n)$ Kac-Moody superalgebra and the definition of affine Lie superalgebra
$A^{(1)}(m-1|n-1)=sl^{(1)}(m|n)$

\begin{definition}

 Suppose that $g$ is $sl(m|n)$. Then, we set a Lie superalgebra $\tilde{g}$ as $g\otimes \mathbb{C}[t^{\pm 1}]\oplus\mathbb{C}c\oplus\mathbb{C}d$ whose commutator relations are following;
\[
[a\otimes t^s,b\otimes t^u]=[a,b]\otimes t^{s+u}+\delta_{s+u}\kappa(a,b)c,
\]
\begin{center}
{c is a central element of $\tilde{g}$}
\end{center}
\[
[d,a\otimes t^s]=sa\otimes t^s
\]
\end{definition}
We also set a superalgebra $\hat{g} \subset \tilde{g}$ as $g\otimes\mathbb{C}[t^{\pm 1}]\oplus\mathbb{C}c$
We have following presentation of $\hat{sl}(m|n)$ as Kac-Moody Lie superalgebra.
Lie Superalgebra $\tilde{sl}(m|n)$ isomorphic to Lie superalgebra $sl^{(1)}(m|n)$ defined by generators $\{x_i^{\pm},h_i,d|0\leq i \leq m+n-1\}$ and following relations:
\[
[d,h_i]=0,\quad [d,x_i^+]=\delta_{i0}x_i^+\quad [d,x_i^-]=-\delta_{i0}x_i^-
\]
\begin{equation}
\label{eq:2}
[h_i,h_j]=0, \quad [h_i,x_j^{\pm}]=\pm a_{i,j}x_j^{\pm}, \quad [x_i^+,x_j^-]=\delta_{i,j}h_i, \quad ad(x_i^{\pm})^{1+|a_{i,j}|}x_j^{\pm}=0
\end{equation}
\begin{equation}
\label{eq:3}
[x_0^{\pm},x_0^{\pm}]=0, \quad [x_m^{\pm},x_m^{\pm}]=0,
\end{equation}
\begin{equation}
\label{eq:4}
[[x_{m-1}^{\pm},x_{m}^{\pm}],[x_{m}^{\pm},x_{m+1}^{\pm}]] = 0,\quad [[x_{m+n-1}^{\pm},x_{0}^{\pm}], [x_{0}^{\pm},x_{1}^{\pm}]] = 0.
\end{equation}
Where the generators $x_m^{\pm}$ and $x_0^{\pm}$ are odd and all other generators are even.
Here
\begin{equation}
\delta_{ij} = \begin{cases}
    1, & if \quad i=j,\\
    0, & if \quad i \neq j.
    \end{cases}
\end{equation}

Note that Lie superalgebra $\hat{sl}(m|n)$ isomorphic to Lie superalgebra over $\mathbb{C}$ defined by generators $\{x_i^{\pm},h_i,d|0\leq i \leq m+n-1\}$ and relations (\ref{eq:2}-\ref{eq:4})\\
Suppose $\{\alpha_i\}_{0\leq i \leq m+n-1}$ is a set of simple roots of $\tilde{sl}(m|n)$ and $\delta$ and $\theta$ are positive roots such that $\delta = \sum\limits_{0 \leq i \leq m+n-1} \alpha_i$, $\theta = \sum\limits_{1 \leq i \leq m+n-1} \alpha_i$. Hence $\alpha_{0}=\delta-\theta$.
Moreover, we set $\Delta$(resp. $\Delta_+$) as a set of roots(resp. positive roots) of $\tilde{sl}(m|n)$.\\
Consider the system of weights introduced by formulas (\ref{rel:4}). We added root $\alpha_{0}$ to this set and introduced root $\delta$ that is dual to $d$ i.e.: $\langle \delta , d \rangle =1$. $\delta$ byorthogonal to Cartan subalgebra $sl(m|n)$. We can consider the simple root $\alpha_{i},i=1,\ldots,m+n-1$ as difference of adjacent weights relative to a given order. Zero root is given by the formula $\alpha_{0}=\delta - \theta$ and we fix this root. Other roots forms a system of simple roots, that called distinguished system of simple of roots if order defined as above and transforms to another system of simple roots under change of order of weights. We define (symmetric) Cartan matrix as symmetric matrix with elements defined as follows $a_{i,i}=0$ if root is odd, other matrix elements defined by following formula:
\begin{equation}
    a_{ij}=-2\frac{(\alpha_i,\alpha_j)}{(\alpha_i,\alpha_i)}
\end{equation}
We suppose, that, as above, $(\varepsilon_{i},\varepsilon_{j})=\delta_{ij}$, $(\delta_{i},\delta_{j})=\delta_{ij}$, where $\delta_{ij}$ is Kronecker delta and $(\varepsilon_{i},\delta_{j})=0$. Every simple root as mentioned above is difference between adjacent weights relative to a given order. Thus, permutations of weights induces the transform of system of simple roots and system of simple roots defined by order of weights. Note that Cartan matrix of Lie superalgebra $sl(\Pi)$ and it's affine analogue consists of diagonal blocks:
\begin{equation}
   \begin{pmatrix}
     \pm 2 & \mp 1\\
    \mp 1 & \ldots
    \end{pmatrix}
\end{equation} for even roots,

\begin{equation}
   \begin{pmatrix}
      0 &  1\\
    1 & \ldots
    \end{pmatrix}
\end{equation} for odd roots.

\section{Affine super Yangian}

Let $S_n$ be the symmetric group and $\hat{sl}(m|n)$ algebra is associated with indecomposable Cartan matrix $(a_{ij})_{i,j \in I}$ where $I$ is the set of vertices of the Dynkin diagram corresponding to $\hat{sl}(m|n)$. We set $\{a,b\}$ as $ab + (-1)^{p(a)p(b)}ba$. Note that elements of $I$ indexed by numbers $\{0,1,\ldots,m+n-1\}$, we will identify this two sets when we want to fix the order on the set of simple roots. Then following \cite{Ueda} we can define Yangian of $\hat{sl}(m|n)$ in distinguished realization.
\begin{definition}
The Yangian $Y_\hbar(\hat{sl}(m|n))$ is unital associative $\mathbb{C}[h]$-algebra generated by the elements $x_{\alpha_{i},r}^{\pm}$, $h_{\alpha_{i},r}$, for $i \in \{1,\ldots,m+n-1 \}$ and $r \in \mathbb{Z}_{\geq 0}$, subject to the relations
\begin{equation}
\label{eq:8}
    [h_{{i},r},h_{{j},s}]=0,
\end{equation}
\begin{equation}
    \label{rel:1}
    [h_{{i},0},x_{{j},s}^{\pm}]=\pm a_{ij}x_{{j},s}^{\pm},
\end{equation}
\begin{equation}
\label{rel:2}
    [x_{{i},r}^{+},x_{{j},s}^{-}]=\delta_{ij}h_{{i},r+s},
\end{equation}
\begin{equation}
    [h_{{i},r+1},x_{{j},s}^{\pm}]-[h_{{i},r},x_{{j},s+1}^{\pm}]=\pm\frac{\hbar a_{ij}}{2}\{h_{{i},r},x_{{j},s}^{\pm}\},
\end{equation}
\begin{equation}
    [x_{{i},r+1}^{\pm},x_{{j},s}^{\pm}]-[x_{{i},r}^{\pm},x_{{j},s+1}^{\pm}]=\pm \frac{ \hbar a_{ij}}{2}\{x_{{i},r}^{\pm},x_{{j},s}^{\pm}\},
\end{equation}
\begin{equation}
    \sum\limits_{\sigma_i \in S_n}[x_{{i},r_{\sigma(i)}}^{\pm},[x_{{i},r_{\sigma(2)}}^{\pm},\ldots,[x_{{i},r_{\sigma(m)}}^{\pm},x_{{j},s}^{\pm}]\ldots]]=0 \quad for \quad i\neq j \quad and \quad n=1-a_{ij},
\end{equation}
\begin{equation}
\label{equ:17}
  [x_{{i},r},x_{{i},s}]=0 \quad (i=0,m),
\end{equation}
\begin{equation}
\label{eq:9}
    [[x_{{i-1},0}^{\pm},x_{{i},0}^{\pm}],[x_{{i},0}^{\pm},x_{{i+1},0}^{\pm}]]=0 \quad (i=0,m)
\end{equation}
where
we set $p:\{1,\ldots ,m+n\} \rightarrow \{0,1\}$ as
\[
p(i) =
    \begin{cases}
      0 & (1\leq i\leq m), \\
      1 & (m+1\leq i\leq m+n).\\
    \end{cases}
\]
\[
a_{i,j} =
    \begin{cases}
      (-1)^{p(i)}+(-1)^{p(i+1)}, & if \quad i=j, \\
      -(-1)^{p(i+1)}, & if \quad j=i+1,\\
      -(-1)^{p(i)}, & if \quad j=i-1,\\
      1, & if \quad (i,j) = (0,m+n-1),(m+n-1,0),\\
      0, & otherwise.\\
    \end{cases}
\]
and the generators $x_{{m},r}^{\pm}$ and $x_{{0},r}^{\pm}$ are odd and all over generators are even.
\end{definition}

We can define affine super Yangian by following equivalent way.
\begin{definition}
Super Yangian $Y_\hbar(\hat{sl}(m|n))$ is unital associative $\mathbb{C}[h]$-algebra generated by the elements $x_{\alpha_{i},r}^{\pm}$, $h_{\alpha_{i},r}$, for $i \in \{1,\ldots,m+n-1 \}$ and $r \in \mathbb{Z}_{\geq 0}$, subject to the relations (\ref{eq:8}-\ref{equ:17}) and the following relation equivalent to (\ref{eq:9})
\begin{equation}
     [[x_{{i-1},k}^{\pm},x_{{i},0}^{\pm}],[x_{{i},0}^{\pm},x_{{i+1},t}^{\pm}]]=0 \quad (i=0,m).
\end{equation}
\end{definition}
Note that here and above $a_{ij}$  are elements of Cartan matrix $A={(a_{ij})}_{i,j \in I}$ of superalgebra $\hat{sl}(m|n)$ defined by distinguished system of simple roots.\\
We also define affine super Yangian of superalgebra $\hat{sl(m|n)}=A^{(1)}(m-1|n-1)$.
\begin{definition}
Suppose that $m,n\geq 2$ and $m\neq n$. We define $Y_{\hbar}(\tilde{sl}(m|n))$ is the associative superalgebra over $\mathbb{C}$ generated by $\{x_{\alpha_i,r}^{\pm}, h_{\alpha_i,r}, d|i\in \{0,1,\ldots,m+n-1\},r\in \mathbb{Z}_{\geq 0}\}$ with parameter $\hbar \in \mathbb{C}$ subject to the relations (\ref{eq:8}-\ref{eq:9}) and

\begin{equation}
    [d, h_{i,r}]=0,\quad  [d, x_{i,r}^{\pm}]= \pm \delta_{i0}x_{i,r}^{\pm},
\end{equation}

where the generators $x_{m,r}^{\pm}$ and $x_{0,r}^{\pm}$ are odd and all over generators are even.\\
\end{definition}
Let us set $\tilde{h}_{\alpha_i,1} = h_{\alpha_i,1}-\frac{\hbar}{2}h_{\alpha_i,0}^2$, then
we can show (\cite{Ueda}) that following algebra is isomorphic to $Y_{\hbar}(\hat{sl}(m|n))$.

\begin{theorem}\label{thm_3.1}

Suppose $m,n\geq 2$ and $m\neq n$. The affine super Yangian $Y_{\hbar}(\hat{sl}(m|n))$ is isomorphic to associative superalgebra generated by $x_{{i},r}^{\pm}$, $h_{i,r}$, for $i \in \{1,\ldots,m+n-1 \}$ and $r \in \{0,1 \}$, subject to the relations:
\begin{equation}
    [h_{{i},r},h_{{j},s}]=0,
\end{equation}
\begin{equation}
    [x_{{i},0}^{+},x_{{j},0}^{-}]=\delta_{ij}h_{i,0},
\end{equation}
\begin{equation}
    [x_{{i},1}^{+},x_{{j},0}^{-}]=\delta_{ij}h_{{i},1}=[x_{{i},0}^{+},x_{{j},1}^{-}],
\end{equation}
\begin{equation}
    [h_{{i},0},x_{i,r}^{\pm}]=\pm a_{ij}x_{j,r}^{\pm},
\end{equation}
\begin{equation}
    [x_{{i},1}^{\pm},x_{{j},0}^{\pm}]-[x_{{i},0}^{\pm},x_{{j},1}^{\pm}]=\pm \frac{ \hbar a_{ij}}{2}\{x_{{i},0}^{\pm},x_{{j},0}^{\pm}\},
\end{equation}
\begin{equation}
    [\tilde{h}_{i,1},x_{j,0}^{\pm}]=\pm a_{ij}x_{j,1}^{\pm},
\end{equation}
\begin{equation}
    (\text{ad} x_{i,0}^{\pm})^{(1+|a_{ij}|)}(x_{j,0}^{\pm}) = 0, \quad (i \neq j),
\end{equation}
\begin{equation}
    [x_{i,0}^{\pm},x_{i,0}^{\pm}]=0, \quad if \quad i=(0,m),
\end{equation}
\begin{equation}
    [[x_{{i-1},0}^{\pm},x_{i,0}^{\pm}][x_{i,0}^{\pm},x_{{i+1},0}^{\pm}]]=0, \quad if \quad i=(0,m).
\end{equation}
\end{theorem}
Suppose $\Pi$ is an arbitrary system of roots of special linear Kac-Moody superalgebra. We define affine super Yangian for arbitrary realization of affine special linear Kac-Moody superalgebra $\hat{sl}(E,\Pi,p)$, where $E$ is vector space, $\Pi=\{\alpha_0,\alpha_1,\ldots,\alpha_{m+n-1}\}$ is a basis of simple roots and $p$ is the parity function defined on root lattice. We define affine super Yangian as above with Cartan matrix $A={(a_{ij})}_{i,j \in I}$ which is defined by the arbitrary system of simple roots $\Pi$.
\begin{definition}
The Yangian $Y_\hbar(\hat{sl}(E,\Pi,p))$ is a unital associative $\mathbb{C}[h]$-algebra generated by the elements $x_{\alpha_{i},r}^{\pm}$, $h_{\alpha_{i},r}$, for $i \in \{1,\ldots,m+n-1 \}$ and $r \in \mathbb{Z}_{\geq 0}$, subject to the relations
\begin{equation}
\label{rel:27}
    [h_{\alpha_{i},r}, h_{\alpha_{j},s}] = 0,
\end{equation}
\begin{equation}
\label{rel:28}
    [h_{\alpha_{i},0},x_{\alpha_{j},s}^{\pm}]=\pm a_{ij}x_{\alpha_{j},s}^{\pm},
\end{equation}
\begin{equation}
\label{rel:29}
    [x_{\alpha_{i},r}^{+},x_{\alpha_{j},s}^{-}]=\delta_{ij}h_{\alpha{i},r+s},
\end{equation}
\begin{equation}
\label{rel:30}
    [h_{\alpha_{i},r+1},x_{\alpha_{j},s}^{\pm}]-[h_{\alpha_{i},r},x_{\alpha_{j},s+1}^{\pm}]=\pm\frac{\hbar a_{ij}}{2}\{h_{\alpha_{i},r},x_{\alpha_{j},s}^{\pm}\},
\end{equation}
\begin{equation}
\label{rel:31}
    [x_{\alpha_{i},r+1}^{\pm},x_{\alpha_{j},s}^{\pm}]-[x_{\alpha_{i},r}^{\pm},x_{\alpha_{j},s+1}^{\pm}]=\pm \frac{ \hbar a_{ij}}{2}\{x_{\alpha_{i},r}^{\pm},x_{\alpha_{j},s}^{\pm}\},
\end{equation}
\begin{equation}
\label{rel:32}
\sum\limits_{\sigma_i \in S_n}[x_{\alpha_{i},r_{\sigma(i)}}^{\pm},[x_{\alpha_{i},r_{\sigma(2)}}^{\pm},\ldots,[x_{\alpha_{i},r_{\sigma(m)}}^{\pm},x_{\alpha_{j},s}^{\pm}]\ldots]]=0, \quad for \quad i\neq j \quad and \quad n=1-a_{ij},
\end{equation}
\begin{equation}
\label{rel:33}
  [x_{\alpha_{i},r},x_{\alpha_{i},s}]=0, \quad \text{for every odd root $\alpha_i$},
\end{equation}
\begin{equation}
\label{rel:34}
[[x_{\alpha_{i-1},r}^{\pm},x_{\alpha_{i},0}^{\pm}],[x_{\alpha_{i},0}^{\pm},x_{\alpha_{i+1},s}^{\pm}]]=0, \quad \text{for every odd root $\alpha_i$}.
\end{equation}
\end{definition}
Here $a_{ij}$ are elements of Cartan matrix defined by given system of simple roots $\Pi$.

\begin{definition}
Suppose $\tilde{sl}(E,\Pi,p)$ is the $A(m|n)$ type of algebra with $m\neq n$ and $m,n\geq 2$. We define $Y_{\hbar}(\tilde{sl}(E,\Pi,p))$ as the associative superalgebra over $\mathbb{C}$ generated by $\{h_{\alpha_i, r},   x_{\alpha_i,r}^{\pm}, d |r\in \mathbb{Z}_{\geq 0}\}$ with parameter $\hbar \in \mathbb{C}$ subject to the relations (\ref{rel:27}-\ref{rel:34}) and

\begin{equation}\label{in:1}
    [d,h_{\alpha_i,r}]=0,\quad  [d, x_{\alpha_i,r}^{\pm}]= \pm \delta_{i0}x_{\alpha_i,r}^{\pm},
\end{equation}

where the odd generators $x_{\alpha_{i},r}^{\pm}$ corresponds odd roots $\alpha_i$ and all other generators are even.\\
\end{definition}




\subsection{Main theorem}

\begin{theorem} \label{thm_3.2}
Suppose $\tilde{sl}(E,\Pi,p)$ is the Lie superalgebra of $A(m|n)$ type with $m\neq n$ and $m,n\geq 2$. The affine super Yangian $Y_{\hbar}(\hat{sl}(E,\Pi,p))$ is isomorphic to associative superalgebra generated by $x_{\alpha_{i},r}^{\pm}$, $h_{\alpha_{i},r}$, for $i \in \{1,\ldots,m+n-1 \}$ and $r \in \{0,1 \}$, subject to the relations:
\begin{equation}
\label{rel:35}
    [h_{\alpha_{i},r},h_{\alpha_{j},s}]=0,
\end{equation}
\begin{equation}
\label{rel:36}
    [x_{\alpha_{i},0}^{+},x_{\alpha_{j},0}^{-}]=\delta_{ij}h_{\alpha_i,0},
\end{equation}
\begin{equation}
\label{rel:37}
    [x_{\alpha_{i},1}^{+},x_{\alpha_{j},0}^{-}]=\delta_{ij}h_{\alpha_{i},1}=[x_{\alpha_{i},0}^{+},x_{\alpha_{j},1}^{-}],
\end{equation}
\begin{equation}
\label{rel:38}
    [h_{\alpha_{i},0},x_{\alpha_j,r}^{\pm}]=\pm a_{ij}x_{\alpha_j,r}^{\pm},
\end{equation}
\begin{equation}
\label{rel:39}
    [x_{\alpha_{i},1}^{\pm},x_{\alpha_{j},0}^{\pm}]-[x_{\alpha_{i},0}^{\pm},x_{\alpha_{j},1}^{\pm}]=\pm \frac{ \hbar a_{ij}}{2}\{x_{\alpha_{i},0}^{\pm},x_{\alpha_{j},0}^{\pm}\},
\end{equation}
\begin{equation}
\label{rel:40}
    [\tilde{h}_{\alpha_i,1},x_{\alpha_j,0}^{\pm}]=\pm a_{ij}x_{\alpha_j,1}^{\pm},
\end{equation}
\begin{equation}
\label{rel:41}
    (\text{ad} x_{\alpha_i,0}^{\pm})^{(1+|a_{ij}|)}(x_{\alpha_j,0}^{\pm}) = 0\quad (i \neq j),
\end{equation}
\begin{equation}
\label{rel:42}
    [x_{\alpha_i,0}^{\pm},x_{\alpha_i,0}^{\pm}]=0 ],\quad \text{for every odd root $\alpha_i$},
\end{equation}
\begin{equation}
\label{rel:43}
    [[x_{\alpha_{i-1},0}^{\pm},x_{\alpha_i,0}^{\pm}][x_{\alpha_i,0}^{\pm},x_{\alpha_{i+1},0}^{\pm}]]=0, \quad \text{for every odd root $\alpha_i$}.
\end{equation}
\end{theorem}

\subsection{Weyl groupoid}

Let $s$ be the element of Weyl groupoid $\hat{W}$. $\hat{W}$ has natural action on system of simple roots $\Pi$ of Lie super algebra $A(m|n)$, for $s \in \hat{W}$ $s:\Pi \rightarrow \Pi_1$. Thus we can define action of Weyl groupoid $\hat{W}$ on super Yangian $\hat{W}:Y(\hat{sl}(E,\Pi,p))\rightarrow Y(\hat{sl}(E,\Pi_1,p))$. Elements $w$ of Weyl groupoid define isomorphisms $T_{w}$ in category of super Yangians $Y(\hat{sl}(E,\Pi,p))$, that we define as quantum reflections. Thus $\hat{W}$ is a groupoid in category sense.

Note that every even element $s\in \hat{W}$ defines automorphism $T_s:Y(\hat{sl}(E,\Pi,p))\rightarrow Y(\hat{sl}(E,\Pi_1,p)$.
Even reflections form Weyl group, elements of Weyl group induces automorphisms $T_{s}$ of super Yangian $Y(\hat{sl}(E,\Pi,p))$.We obtained naturally defined action of Weyl group on any super Yangian, that we observe as the element of category noticed above.

Note that this definition of Weyl groupoid is  same for case of affine super Yangian $Y(\hat{sl}(E,\Pi,p))$.

\begin{theorem}
For every element $s \in \hat{W}$ exists isomorphism $T_s:Y(\hat{sl}(E,\Pi,p))\rightarrow Y(\hat{sl}(E,\Pi_1,p)$ such that $T_s$ is an automorphism if and only if $s$ is an even reflection.
\end{theorem}

For even simple reflections we can define map $T_{\alpha_j}:Y(\hat{sl}(E,\Pi,p))\rightarrow Y(\hat{sl}(E,\Pi,p))$ as done in Kodera work \cite{Kodera}
\[
T_{\alpha_i}(x_{\alpha_j,0}^{\pm})=
\begin{cases}
         -x_{\alpha_i,0}^{\mp}, & if \quad a_{i,j}=2,\\
        \pm[x_{\alpha_i,0}^{\pm},x_{\alpha_j,0}^{\pm}], & if \quad a_{i,j}=-1,\\
         x_{\alpha_j,0}^{\pm}, & if \quad a_{i,j}=0,\\
\end{cases}
\]
\[
T_{\alpha_i}(h_{\alpha_j,0})=
\begin{cases}
    -h_{\alpha_i,0}, & if \quad a_{i,j}=2,\\
    h_{\alpha_i,0}+h_{\alpha_j,0}, & if \quad a_{i,j}=-1, \\
    h_{\alpha_j,0}, & if \quad a_{i,j}=0.
\end{cases}
\]
\[
T_{\alpha_i}(x_{\alpha_j,1}^{\pm})=
\begin{cases}
    -x_{\alpha_i,1}^{\mp}+\frac{\hbar}{2}\{h_{\alpha_i,0},x_{\alpha_i,0}^{\mp}\}, & if \quad a_{i,j}=2,\\
    \pm[x_{\alpha_i,0}^{\pm},x_{\alpha_j,1}^{\pm}], & if \quad a_{i,j}=-1,\\
     x_{\alpha_j,1}^{\pm}, & if \quad a_{i,j}=0,\\
\end{cases}
\]
\[
T_{\alpha_i}(\tilde{h}_{\alpha_j,1})=
\begin{cases}
    -\tilde{h}_{\alpha_j,1}-\hbar\{x_{\alpha_i,0}^{+},x_{\alpha_i,0}^{-}\},& if \quad a_{i,j}=2,\\
    \tilde{h}_{\alpha_j,1}+\tilde{h}_{\alpha_i,1}+\frac{\hbar}{2}\{x_{\alpha_i,0}^{+},x_{\alpha_i,0}^{-}\},& if \quad a_{i,j}=-1, \\
    \tilde{h}_{\alpha_j,1},& if \quad a_{i,j}=0.\\
\end{cases}
\]
Let us define quantum odd reflections. Let $\alpha_{i}$ simple odd root. Let $s_{\alpha_{i}}:\Pi \rightarrow \Pi_{1}$ be the reflection induced by this root(odd reflection). Let $\beta_{j}:=s_{\alpha_{i}}(\alpha_{j})\in \Pi_1$ be the image of simple root $\alpha_{j}$ of system of simple roots $\Pi$ under action of reflection $s_{\alpha_{i}}=s_{i}$.

For simple odd reflections we define map $T_{\alpha_j}:Y(\hat{sl}(E,\Pi,p))\rightarrow Y(\hat{sl}(E,\Pi,p))$ as follows

\[
T_{\alpha_i}(x_{\alpha_j,0}^{\pm})=
\begin{cases}
         -x_{s_i(\alpha_i),0}^{\mp}, & if \quad i=j,\\
        \pm[x_{s_i(\alpha_i),0}^{\pm},x_{s_i(\alpha_j),0}^{\pm}], & if \quad a_{i,j}=-1,\\
         x_{s_i(\alpha_j),0}^{\pm}, & if \quad a_{i,j}=0\quad i\neq j,\\
\end{cases}
\]
\[
T_{\alpha_i}(h_{\alpha_j,0})=
\begin{cases}
    -h_{s_i(\alpha_i),0}, & if \quad i=j,\\
    h_{s_i(\alpha_i),0}+h_{s_i(\alpha_j),0}, & if \quad a_{i,j}=-1, \\
    h_{s_i(\alpha_j),0}, & if \quad a_{i,j}=0, \quad i \neq j,
\end{cases}
\]

\[
T_{\alpha_i}(x_{\alpha_j,1}^{\pm})=
\begin{cases}
    -x_{s_i(\alpha_i),1}^{\mp}, & if \quad i=j,\\
    \pm[x_{s_i(\alpha_i),0}^{\pm},x_{s_i(\alpha_j),1}^{\pm}], & if \quad a_{i,j}=-1,\\
     x_{s_i(\alpha_j),1}^{\pm}, & if \quad a_{i,j}=0, \quad i \neq j,\\
\end{cases}
\]
\[
T_{\alpha_i}(\tilde{h}_{\alpha_j,1})=
\begin{cases}
    -\tilde{h}_{s_i(\alpha_j),1}, & if \quad i=j,\\
    \tilde{h}_{s_i(\alpha_j),1}+\tilde{h}_{s_i(\alpha_i),1}+\frac{\hbar}{2}\{x_{s_{i}(\alpha_{i}),0}^{+},x_{s_{i}(\alpha_{i}),0}^{-}\}, &if \quad a_{i,j}=-1, \\
    \tilde{h}_{s_i(\alpha_j),1}, & if \quad a_{i,j}=0, \quad i \neq j.\\
\end{cases}
\]

\section{Proofs}

\subsection{Proof of theorem  \ref{thm_3.1}}

Super algebra introduced in theorem \ref{thm_3.2} we denote as $Y_{\hbar}^1(\hat{sl}(E,\Pi,p))$. By the definition of $\tilde{h}_{\alpha_i,1}$ we can write (\ref{rel:30}) as
\begin{equation}
\label{rel:44}
    [\tilde{h}_{\alpha_i,1},x_{\alpha_j,r}^{\pm}]=\pm a_{ij}x_{\alpha_j,r+1}^{\pm}
\end{equation}
By relations (\ref{rel:44}) and (\ref{rel:29}) we can write following relations:
\begin{equation}
    x_{\alpha_i,r+1}^{\pm}=\pm \frac{1}{a_{ii}}[\tilde{h}_{\alpha_{i+1},1},x_{\alpha_i,r}], \quad h_{\alpha_1,r+1}=[x_{\alpha_i,r+1}^{+},x_{\alpha_i,0}^{-}], \quad \text{if $\alpha_i$ is even root},
\end{equation}
\begin{equation}
\label{rel:46}
    x_{\alpha_i,r+1}^{\pm}=\pm \frac{1}{a_{i+1,i}}[\tilde{h}_{\alpha_{i+1},1},x_{\alpha_i,r}], \quad h_{\alpha_1,r+1}=[x_{\alpha_i,r+1}^{+},x_{\alpha_i,0}^{-}], \quad \text{if $\alpha_i$ is odd root}
\end{equation}
for all $r \geq 1$.

\begin{lemma} \label{lem:1}
The relation (\ref{rel:1}) holds for all $i,j\in I$ in $Y_{\hbar}^1(\hat{sl}(E,\Pi,p))$. For all $i,j\in I$ we obtain
\begin{equation}
\label{rel:47}
    [\tilde{h}_{\alpha_i,1},x_{\alpha_j,r}^{\pm}]=\pm a_{ij}x_{\alpha_j,r+1}^{\pm}
\end{equation}
in $Y_{\hbar}^1(\hat{sl}(E,\Pi,p))$\\
\end{lemma}
We show only the case then $\alpha_j$ is the odd root. We prove this lemma by induction on $r$. For $r=0$, the relations is just relations (\ref{rel:38}), (\ref{rel:40}).Suppose that relations (\ref{rel:1}), (\ref{rel:44}) holds for $r=k$. First, we show that relation \ref{rel:1} holds than $r=k+1$. By (\ref{rel:44}) we obtain
\begin{equation}
\label{rel:48}
    [h_{\alpha_i,0},x_{\alpha_j,k+1}^{\pm}]=\pm \frac{1}{a_{j,j+1}}[h_{\alpha_i,0}[\tilde{h}_{\alpha_{j+1},1},x_{\alpha_j,k}]].
\end{equation}

By $[h_{\alpha_i,0},h_{\alpha_j,1}]=0$, we find that right hand sight of (\ref{rel:48}) is
\begin{equation}
    \pm \frac{1}{a_{j,j+1}}[\tilde{h}_{\alpha_{j+1},1},[h_{\alpha_i,0},x_{\alpha_j,k}^{\pm}]]=\frac{a_{i,j}}{a_{j,j+1}}[\tilde{h}_{\alpha_{j+1},1},x_{\alpha_j,k}^{\pm}]=\frac{a_{i,j}}{a_{j,j+1}}(\pm a_{j,j+1}x_{\alpha_j,k+1}^{\pm})=a_{i,j}x_{\alpha_j,k+1}^{\pm}.
\end{equation}
Thus we shown that $[h_{\alpha_i,0},x_{\alpha_j,k+1}]=\pm a_{ij}x_{\alpha_j,k+1}$.
Next we show that (\ref{rel:1}). We check the relation
\begin{equation}
    [h_{\alpha_i,1},x_{\alpha_j,k+1}^{\pm}]=\pm a_{ij}x_{\alpha_j,k+2}.
\end{equation}
By (\ref{rel:46}) we obtain
\begin{equation}
\label{rel:51}
     [h_{\alpha_i,1},x_{\alpha_j,k+1}^{\pm}]=\pm \frac{1}{a_{j,j+1}}[\tilde{h}_{\alpha_i,1},[\tilde{h}_{\alpha_{j+1},1},x_{\alpha_j,k}]].
\end{equation}
By $[h_{\alpha_i,1},h_{\alpha_j,1}]=0$ we find that right hand side is equal to
\[
\frac{1}{a_{j,j+1}}[\tilde{h}_{\alpha_{j+1},1},[\tilde{h}_{\alpha_i,1},x_{\alpha_j,k}]].
\]
By the induction hypothesis we can write hand side (\ref{rel:51}) as
\begin{equation}
    \label{rel:52}
    \pm \frac{1}{a_{j,j+1}}[\tilde{h}_{\alpha_{j+1},1},[\tilde{h}_{\alpha_i,1},x_{\alpha_j,k}]]=\frac{a_{ij}}{a_{j,j+1}}[\tilde{h}_{\alpha_{j+1},1},x_{\alpha_j,k+1}^{\pm}].
\end{equation}
The generator $x_{\alpha_j,k+2}^{\pm}$ was defined in (\ref{rel:46}), we find that right hand side of (\ref{rel:52}) is $\pm a_{ij}x_{\alpha_j,k+2}^{\pm}$. This statement comletes the proof of lemma 4.1.1.

\begin{lemma}
\label{lem:2}
\begin{enumerate}
    \item The relation (\ref{rel:29}) holds in $Y_{\hbar}^1(\hat{sl}(E,\Pi,p))$. For all $i,j\in I$ when $i=j$ and $r+s \leq 2$.
    \item Suppose that $i\neq j$. Then the relations (\ref{rel:29}), (\ref{rel:31}) hold for any $r,s$ in $Y_{\hbar}^1(\hat{sl}(E,\Pi,p))$.
    \item The relation (\ref{rel:31}) holds in $Y_{\hbar}^1(\hat{sl}(E,\Pi,p))$ when $i=j,(r,s)=(1,0)$.
    \item The relation (\ref{rel:30}) holds in $Y_{\hbar}^1(\hat{sl}(E,\Pi,p))$ when $i=j,(r,s)=(1,0)$.
    \item For all $i,j$, the relation (\ref{rel:30}) holds in $Y_{\hbar}^1(\hat{sl}(E,\Pi,p))$ when $(r,s)=1,0$.
    \item Set $\tilde{h}_{\alpha_i,2}=h_{\alpha_i,2}-h_{\alpha_i,0}h_{\alpha_i,1}+\frac{1}{3}h_{\alpha_i,0}^3$. Then the following equation holds for all $i=j$ in $Y_{\hbar}^1(\hat{sl}(E,\Pi,p))$;\\
    $[\tilde{h}_{\alpha_i,2},x_{\alpha_j,0}^{\pm}]=\pm a_{ij}x_{\alpha_j,2}^{\pm} \pm \frac{1}{12}  a_{ij}^3x_{\alpha_j,0}^{\pm}$
    \item For all $i,j$ the relation (\ref{rel:32}) holds in $Y_{\hbar}^1(\hat{sl}(E,\Pi,p))$ when
    \begin{itemize}
        \item $r_1=\ldots=r_b=0, s\in \mathbb{Z}_{\geq 0}$,
        \item $r_1=1,r_2=\ldots=r_b=0, s\in \mathbb{Z}_{\geq 0}$,
        \item $r_1=2,r_2=\ldots=r_b=0, s\in \mathbb{Z}_{\geq 0}$,
        \item $b \geq 2 $ and $r_1=r_2,r_3=\ldots=r_b=0, s\in \mathbb{Z}_{\geq 0}$.
    \end{itemize}
    \item In $Y_{\hbar}^1(\hat{sl}(E,\Pi,p))$, we have $[h_{\alpha_j,1},x_{\alpha_i,1}^{\pm}]=\frac{a_{ij}}{a_{ii}}[h_{\alpha_i,1},x_{\alpha_i,1}^{\pm}]\pm \frac{a_{ij}}{2}(\{h_{\alpha_j,0},x_{\alpha_i,1}\}-\{h_{\alpha_i,0},x_{\alpha_i,1} \})$ for all $i,j$ for all even roots.
    \item For all $i,j$ we have $[h_{\alpha_i,2},h_{\alpha_j,0}]$ in $Y_{\hbar}^1(\hat{sl}(E,\Pi,p))$.
    \item Let $a_{ii}=2$ and $a_{ij}=-1$, Then $[h_{\alpha_i,2},h_{\alpha_i,1}]$ holds in $Y_{\hbar}^1(\hat{sl}(E,\Pi,p))$.
\end{enumerate}
\end{lemma}

We prove statements (1-5), since the proof of statements (6-10) contained in the set of lemmas(2.33-2.36) in \cite{Guaywork}. the proofs of of first and second statement are the same as those of Lemma 2.22 and 2.26 in \cite{Guaywork} in case where roots are odd the proof of third fourth and fifth statement of Lemma (\ref{lem:2}) are similar to proofs of Lemma 2.23, Lemma 2.24, Lemma 2.28 in \cite{Guaywork}. We omit it.
Let $i=j$ and $\alpha_i$ is an odd root, then we apply $ad(\tilde{h}_{\alpha_{i+1},1})$ to (\ref{rel:42}), then we have $a_{i,i+1}[x_{\alpha_i,1}^{\pm},x_{\alpha_i,0}^{\pm}]\pm a_{i,i+1}[x_{\alpha_i,0}^{\pm},x_{\alpha_i,1}^{\pm}]$. Since $[x_{\alpha_i,0}^{\pm},x_{\alpha_i,1}^{\pm}]=[x_{\alpha_i,1}^{\pm},x_{\alpha_i,0}^{\pm}]$ we obtain $[x_{\alpha_i,0}^{\pm},x_{\alpha_i,1}]=[x_{\alpha_i,1}^{\pm},x_{\alpha_i,0}^{\pm}]=0$. Next, we show that $[x_{\alpha_i,2}^{\pm},x_{\alpha_i,0}^{\pm}]=[x_{\alpha_i,1}^{\pm},x_{\alpha_i,1}^{\pm}]=[x_{\alpha_i,0}^{\pm},x_{\alpha_i,2}^{\pm}]$ holds. We apply $ad(\tilde{h}_{\alpha_{i+1},1})$ to $[x_{\alpha_i,0}^{\pm},x_{\alpha_i,1}^{\pm}]=[x_{\alpha_i,1}^{\pm},x_{\alpha_i,0}^{\pm}]$, we obtain
\begin{equation}
\label{rel:53}
    \pm a_{i,i+1}([x_{\alpha_i,2}^{\pm},x_{\alpha_i,0}^{\pm}]+[x_{\alpha_i,1}^{\pm},x_{\alpha_i,1}^{\pm}])=\pm a_{i,i+1}([x_{\alpha_i,0}^{\pm},x_{\alpha_i,2}^{\pm}]+[x_{\alpha_i,1}^{\pm},x_{\alpha_i,1}^{\pm}])=0.
\end{equation}
Let $\alpha_j$ be an odd root and $\alpha_i=\alpha_{j+1}$, we can prove fifth statement of Lemma (\ref{lem:2}) similar as Lemma 2.28 in \cite{Guaywork}. Then similar to Lemma 1.4 in \cite{Levendorski} there exists $\hat{h}_{\alpha_{i+1},2}$ such that
\begin{equation}
    [\hat{h}_{\alpha_{i+1},2},x_{\alpha_i,0}^{\pm}]=\pm a_{i,i+1}x_{\alpha_i,2}^{\pm}.
\end{equation}
Applying $ad(\hat{h}_{\alpha_{i+1},2})$ to (\ref{rel:42}), we obtain
\begin{equation}
\label{rel:55}
    \pm a_{i,i+1}([x_{\alpha_i,0}^{\pm},x_{\alpha_i,2}^{\pm}]+[x_{\alpha_i,0}^{\pm},x_{\alpha_i,2}^{\pm}]) = 0.
\end{equation}

Equations (\ref{rel:53}), (\ref{rel:55}) are linearly independent, thus we obtain that $[x_{\alpha_i,2}^{\pm},x_{\alpha_i,0}^{\pm}]=[x_{\alpha_i,1}^{\pm},x_{\alpha_i,1}^{\pm}]=[x_{\alpha_i,0}^{\pm},x_{\alpha_i,2}^{\pm}]$, thus we proved statement (3).

\begin{lemma}
\label{lem:3}
Suppose $\alpha_i$ be an even root and, then we have
$[h_{\alpha_i,2},h_{\alpha_i,1}]=0$ in $Y_{\hbar}^1(\hat{sl}(E,\Pi,p))$.\\
\end{lemma}

We change $h_{\alpha_i,r}, x_{\alpha_i,r}^+$ in Proposition 2.36 in \cite{Guaywork} into $-h_{\alpha_i,r}, -x_{\alpha_i,r}^+$.Then we obtain the $[-h_{\alpha_i,2},-h_{\alpha_i,1}]=0$. That completes the proof.\\

We know that $[h_{\alpha_i,2},h_{\alpha_i,1}]=0$ holds for every even root in $I$. By Lemma (\ref{lem:2}) statement 10 and Lemma (\ref{lem:3}) we prove the following lemma in the same way as Proposition 2.39 in \cite{Guaywork}. The condition need to be completed is existence of one even root, since our algebra is affinization of $\hat{sl}(E,\Pi,p)$ with $m,n\geq 2$ and $m\neq n$ this condition is always completed.

\begin{lemma}
\label{lem:4}
Suppose that $\alpha_i$ is odd root and $(\alpha_i,\alpha_j)\neq 0$, then we have $[h_{\alpha_j,2},h_{\alpha_j,1}]=0$ in $Y_{\hbar}^1(\hat{sl}(E,\Pi,p))$.\\
\end{lemma}

By using the relation  $[h_{\alpha_i,2},h_{\alpha_i,1}]=0$ we prove following lemma similar to Theorem 1.2 in \cite{Levendorski}. Since our algebra is the affinization of the type $A(m|n)$ with $m \neq n$ and $m,n\geq 2$ we have at least one even root, hence this condition completes the proof since the only condition is existence of root with $a_{ii}\neq 0$.

\begin{lemma}
\label{lem:5}
\begin{enumerate}
\item The relation (\ref{rel:27}) holds in $Y_{\hbar}^1(\hat{sl}(E,\Pi,p))$ for every even root $\alpha_i$.
\item The relation (\ref{rel:28}) holds in $Y_{\hbar}^1(\hat{sl}(E,\Pi,p))$ for every even root $\alpha_i$.
\item The relation (\ref{rel:31}) holds in $Y_{\hbar}^1(\hat{sl}(E,\Pi,p))$ for every even root $\alpha_i$.
\item The relation (\ref{rel:30}) holds in $Y_{\hbar}^1(\hat{sl}(E,\Pi,p))$ for every even root $\alpha_i$.
\end{enumerate}
\end{lemma}

Next we prove same statement as that of Lemma (\ref{lem:5}) in case then $\alpha_j$ is odd root.
\begin{lemma}
\label{lem:6}
\begin{enumerate}
\item The relation (\ref{rel:31}) holds in $Y_{\hbar}^1(\hat{sl}(E,\Pi,p))$ when $\alpha_j=\alpha_i$ is odd root, hence the relation (\ref{rel:33}) holds in $Y_{\hbar}^1(\hat{sl}(E,\Pi,p))$.
\item The relation (\ref{rel:28}) holds in $Y_{\hbar}^1(\hat{sl}(E,\Pi,p))$ when $\alpha_j=\alpha_i$ is an odd root.
\item We have $[h_{\alpha_i,r},x_{\alpha_i,0}]=0$ when $\alpha_i$ is an odd root in $Y_{\hbar}^1(\hat{sl}(E,\Pi,p))$.
\item The relation (\ref{rel:30}) holds in $Y_{\hbar}^1(\hat{sl}(E,\Pi,p))$ when $\alpha_j=\alpha_i$ is an odd root.
\item The relation (\ref{rel:27}) holds in $Y_{\hbar}^1(\hat{sl}(E,\Pi,p))$ when $\alpha_j=\alpha_i$ is an odd root.
\end{enumerate}
\end{lemma}

We show that $[x_{\alpha_i,r}^{\pm},x_{\alpha_i,s}^{\pm}]=0$ holds. We prove that (\ref{rel:31}) holds by the induction on k=r+s. When $k=0$ this is relation (\ref{rel:42}). We apply $ad(\tilde{h}_{\alpha_{i+1}})$ to $[x_{\alpha_i,0}^+,x_{\alpha_i,0}^+]$, we obtain $a_{i+1,i}([x_{\alpha_i,1}^+,x_{\alpha_i,0}^+]+[x_{\alpha_i,0}^+,[x_{\alpha_i,1}^+])=0$. As we already proved that $[x_{\alpha_i,1}^+,x_{\alpha_i,0}^+]=[x_{\alpha_i,0}^+,[x_{\alpha_i,1}^+]$, we obtain $[x_{\alpha_i,1}^+,x_{\alpha_i,0}^+]=[x_{\alpha_i,0}^+,x_{\alpha_i,1}^+]=0$.
Suppose that $[x_{\alpha_i,r}^{\pm},x_{\alpha_i,s}^{\pm}]=0$ holds for all $r,s$ such that $r+s=k,k+1$. We apply $ad(\tilde{h}_{i+1,i})$ to $[x_{\alpha_i,l}^+,x_{\alpha_i,k+1-l}]=0$, we have
\begin{equation}
\label{rel:56}
    [\tilde{h}_{i+1,i},[x_{\alpha_i,l}^+,x_{\alpha_i,k+1-l}^{+}]]=0.
\end{equation}
Using statement (4) from Lemma (\ref{lem:2}) and induction hypothesis we have
\begin{equation}
\label{rel:57}
[\tilde{h}_{i+1,1},[x_{\alpha_i,l}^+,x_{\alpha_i,k+1-l}]]=a_{i,i+1}([x_{\alpha_i,l+1}^+,x_{\alpha_i,k+1-l}]+[x_{\alpha_i,l}^+,x_{\alpha_i,k+2-l}])
\end{equation}
for all elements of Cartan matrix we have $a_{i,i+1}\neq 0$ since we have algebra $sl(E,\Pi,p)$, thus we have the relation
\begin{equation}
    [x_{\alpha_i,l+1}^+,x_{\alpha_i,k+1-l}]=-[x_{\alpha_i,l}^+,x_{\alpha_i,k+2-l}^{+}]
\end{equation}
by relations (\ref{rel:56}), (\ref{rel:57}), hence we obtain relation
\begin{equation}
\label{rel:59}
    [x_{\alpha_i,l+2}^+,x_{\alpha_i,k-l}^+]=[x_{\alpha_i,l}^+,x_{\alpha_i,k+2-l}^+].
\end{equation}

We apply $ad(\tilde{h}_{\alpha_{i+1},2})$ to $[x_{\alpha_i,l}^+,x_{\alpha_i,k-l}]=0$. Using induction hypothesis and Lemma (\ref{lem:2}) statement (7) and Lemma (\ref{lem:4}) we obtain
\begin{equation}
    [\tilde{h}_{\alpha_{i+1},2},[x_{\alpha_i,l}^+,x_{\alpha_i,k-l}]]=a_{i,i+1}([x_{\alpha_i,l+2}^+,x_{\alpha_i,k-l}^+]+[x_{\alpha_i,l}^+,x_{\alpha_i,k+2-l}^+]).
\end{equation}
Since $a_{i,i+1}\neq 0$, thus we obtain relation
\begin{equation}
\label{rel:61}
    [x_{\alpha_i,l+2}^+,x_{\alpha_i,k-l}^+]=-[x_{\alpha_i,l}^+,x_{\alpha_i,k+2-l}^+].
\end{equation}
Equations (\ref{rel:59}) and (\ref{rel:61}) are linearly independent, thus we obtain that $[x_{\alpha_i,l}^{+},x_{\alpha_i,k+2-l}^{+}]=0$ holds.\\

Proof of statement (2)\\
We prove the statement by induction on $r+s=k$. When $k=0$ this statement become (\ref{rel:42}). Suppose $[x_{\alpha_i,r}^+,x_{\alpha_i,s}^+]=h_{\alpha_i,r+s}$ for all $r,s$ such that $r+s\leq k$, then by definition we have
\begin{equation}
[h_{\alpha_i,r+1},x_{\alpha_{i+1},s}^+]-[h_{\alpha_i,r},x_{\alpha_{i+1},s+1}^+]=[[x_{\alpha_{i},r+1}^+,x_{\alpha_i,0}^-],x_{\alpha_{i},s}^+]-[[x_{\alpha_{i},r}^+,x_{\alpha_i,0}^-],x_{\alpha_{i+1},s+1}^+].
\end{equation}
By Jacobi identity and Lemma (\ref{lem:2}) statement (4) we obtain
\begin{equation}
[h_{\alpha_i,r+1},x_{\alpha_{i+1},s}^+]-[h_{\alpha_i,r},x_{\alpha_{i+1},s+1}^+] = [\{[x_{\alpha_i,r+1}^+,x_{\alpha_{i+1},s}^+]-[x_{\alpha_i,r}^+,x_{\alpha_{i+1},s+1}^+]\},x_{\alpha_i,0}^-].
\end{equation}
By Lemma (\ref{lem:2}) statement (4) we have
\begin{equation}
    [h_{\alpha_i,r+1},x_{\alpha_{i+1},s}^+]-[h_{\alpha_i,r},x_{\alpha_{i+1},s+1}^+] = [\pm a_{i,i+1}\frac{\hbar}{2}\{x_{\alpha_i,r}^+,x_{\alpha_{i+1},s}^+\},x_{\alpha_i,0}^-].
\end{equation}
By using Lemma (\ref{lem:2}) statement (4) we have
\begin{equation}
    [h_{\alpha_i,r+1},x_{\alpha_{i+1},s}^+]-[h_{\alpha_i,r},x_{\alpha_{i+1},s+1}^+]= \pm a_{i,i+1}\frac{\hbar}{2}\{h_{\alpha_i,r},x_{\alpha_{i+1},s+1}\}.
\end{equation}
By the assumption $[x_{\alpha_i,r}^+,x_{\alpha_i,s}^-]=h_{\alpha_i,r+s}$ for $r+s\leq k$ we have
\begin{equation}
    [h_{\alpha_i,r+1},x_{\alpha_{i+1},s}^+]-[h_{\alpha_i,r},x_{\alpha_{i+1},s+1}^+] = [[x_{\alpha_{i},r+1}^+,x_{\alpha_i,0}^-],x_{\alpha_{i},s}^+]-[[x_{\alpha_{i},r}^+,x_{\alpha_i,0}^-],x_{\alpha_{i+1},s+1}^+].
\end{equation}
Since $r+1\leq k$ we obtain
\begin{equation}
[h_{\alpha_i,r+1},x_{\alpha_{i+1},s}^-]-[h_{\alpha_i,r},x_{\alpha_{i+1},s+1}^-]= [[x_{\alpha_{i},r+1}^+,x_{\alpha_i,0}^-],x_{\alpha_{i},s}^-]- [[x_{\alpha_{i},r}^+,x_{\alpha_i,0}^-],x_{\alpha_{i+1},s+1}^-].
\end{equation}
By using Lemma (\ref{lem:2}) statement (4)
\begin{equation}
[h_{\alpha_i,r+1},x_{\alpha_{i+1},s}^-]-[h_{\alpha_i,r},x_{\alpha_{i+1},s+1}^-]= [x_{\alpha_i,r}^+,-a_{i,i+1}\frac{\hbar}{2}\{x_{\alpha_i,0}^-,x_{\alpha_{i+1},s}^-\}].
\end{equation}
Thus by Lemma (\ref{lem:2}) we obtain the relation
\begin{equation}
[h_{\alpha_i,r+1},x_{\alpha_{i+1},s}^-]-[h_{\alpha_i,r},x_{\alpha_{i+1},s+1}^-]= -a_{i+1,i}\frac{\hbar}{2}\{h_{\alpha_i,r},x_{\alpha{i+1},s}^-\}.
\end{equation}
By the similar discussion to Lemma 1.4 in \cite{Levendorski} there exists $\tilde{h}_{\alpha_i,k}$ such that
\[
    \tilde{h}_{\alpha_i,k}=h_{\alpha_i,k}+C(h_{\alpha_i,t}) \quad \{0\leq t \leq k-1\},
\]
where $C(h_{\alpha_i,t})$ is polynom.
\[
[\tilde{h}_{\alpha_i,k},x_{\alpha_{i+1},1}^+]=a_{i,i+1}x_{\alpha_{i+1},k+1}^+, \quad [\tilde{h}_{\alpha_i,k},x_{\alpha_{i+1},1}^-]=-a_{i,i+1}x_{\alpha_{i+1},k+1}^-.
\]
By the assumption that  $[x_{\alpha_i,r}^+,x_{\alpha_i,s}^-]=h_{\alpha_i,r+s}$ for $r+s\leq k$ we have
\[
[\tilde{h}_{\alpha_{i+1},1},h_{\alpha_i,s}]= [[\tilde{h}_{\alpha_{i+1},1},x_{\alpha_i,s}^+],x_{\alpha_i,0}^-]+ [x_{\alpha_i,s}^+,[\tilde{h}_{\alpha_{i+1},1},x_{\alpha_i,0}^-]]=0.
\]
for all $s\leq k$. Thus, it is enough to prove that $[\tilde{h}_{\alpha_{i},k},h_{\alpha_{i+1},1}]=0$ holds. We obtain by the definition of $\tilde{h}_{\alpha_{i+1},1}$
\begin{equation}
[\tilde{h}_{\alpha_{i},k},h_{\alpha_{i+1},1}] = [\tilde{h}_{\alpha_{i},k},[x_{\alpha_{i+1},1}^+,x_{\alpha_{i+1},0}^-]]= a_{i,i+1}[x_{\alpha_{i+1},k+1}^+,x_{\alpha_{i+1},0}^-] - a_{i,i+1}[x_{\alpha_{i+1},1}^+,x_{\alpha_{i+1},k}^-].
\end{equation}
By Lemma (\ref{lem:5}) it is equal to zero.
We apply $ad(\tilde{h}_{\alpha_{i+1},1})$ to $[x_{\alpha_i,r}^+,x_{\alpha_i,k-r}^-]=h_{\alpha_i,k}$. We obtain following relation by induction hypothesis
\begin{equation}
\label{rel:71}
    [\tilde{h}_{\alpha_{i+1},1},[x_{\alpha_i,r}^+,x_{\alpha_i,k-r}^-]] = [\tilde{h}_{\alpha_{i+1},1},h_{\alpha_i,k}].
\end{equation}
By Lemma (\ref{lem:2}) we can rewrite (\ref{rel:71}) as
\begin{equation}
a_{i,i+1}([x_{\alpha_i,r+1}^+,x_{\alpha_i,k-r}^-]- [x_{\alpha_i,r}^+,x_{\alpha_i,k-r+1}^-])= [\tilde{h}_{\alpha_{i+1},1},h_{\alpha_i,k}]=0.
\end{equation}
We proved that $[\tilde{h}_{\alpha_{i+1},1},h_{\alpha_i,k}]=0$.
Proof of statement (3)\\
By statement (2) $[h_{\alpha_i,r},x_{\alpha_i,0}^+]=[[x_{\alpha_i,r}^+,x_{\alpha_i,0}^-],x_{\alpha_i,0}^+]$. By statement (1) and the Jacobi identity, we have
\begin{equation}
\label{rel:73}
    [[x_{\alpha_i,r}^+,x_{\alpha_i,0}^-],x_{\alpha_i,0}^+]=[x_{\alpha_i,r}^+,[x_{\alpha_i,0}^-,x_{\alpha_i,0}^+]].
\end{equation}
The right hand side of (\ref{rel:73}) by the Lemma (\ref{lem:1}) is equal to zero, since root is odd and $[x_{\alpha_i,r}^+,[x_{\alpha_i,0}^-,x_{\alpha_i,0}^+]] = [x_{\alpha_i,r}^+,h_{\alpha_i,0}]$.\\

Proof of statement (4)\\
We prove it by induction on $s$. When $s=0$ it is similar to statement (3). Suppose that $[h_{\alpha_i,r},x_{\alpha_i,s}^+]$ holds. We apply $ad(\tilde{h}_{\alpha_{i+1},1})$ to $[h_{\alpha_i,r},x_{\alpha_i,s}^+]=0$ we obtain
\begin{equation}
\label{rel:74}
    [\tilde{h}_{\alpha_{i+1},1},[h_{\alpha_i,r},x_{\alpha_i,s}^+]]=0.
\end{equation}
By the proof of statement (2) we obtain $[\tilde{h}_{\alpha_{i+1},1},h_{\alpha_i,r}]=0$. Thus the left hand side of (\ref{rel:74}) is $[{h}_{\alpha_{i},r},[\tilde{h}_{\alpha_{i+1},1},x_{\alpha_i,s}^+]]$. By Lemma 4.1.2 statement (4) we obtain
\begin{equation}
\label{rel:75}
    [{h}_{\alpha_{i},r},[\tilde{h}_{\alpha_{i+1},1},x_{\alpha_i,s}^+]]=a_{i,i+1}[h_{\alpha_i,r},x_{\alpha_i,s+1}^+].
\end{equation}
By induction hypothesis, the right hand side of is equal to $a_{i,i+1}[h_{\alpha_i,r},x_{\alpha_i,s+1}^+]$. Since $a_{i,i+1}\neq 0$ we obtain $[h_{\alpha_i,r},x_{\alpha_i,s+1}^+]=0$.For $-$ case proof is similar.\\
Proof of statement (5)\\
By statement (2) $[h_{\alpha_i,r},h_{\alpha_i,s}]=[h_{\alpha_i,r},[x_{\alpha_i,s}^+,x_{\alpha_i,0}^-]$. By the Jacobi identity we have
\begin{equation}
\label{rel:76}
[h_{\alpha_i,r},[x_{\alpha_i,s}^+,x_{\alpha_i,0}^-]]= [[h_{\alpha_i,r},x_{\alpha_i,s}^+], x_{\alpha_i,0}^-]+ [x_{\alpha_i,s}^+,[h_{\alpha_i,r},x_{\alpha_i,0}^-]].
\end{equation}
By statement (4) right hand side of (\ref{rel:76}) is equal to zero. Thus we have shown that $[h_{\alpha_i,r},h_{\alpha_i,s}]=0$.\\
We obtain the relation (\ref{rel:31}) by Lemma (\ref{lem:2}) statement (2), Lemma  (\ref{lem:5}) statement (3) and Lemma (\ref{lem:6}) statement (1). Relation (\ref{rel:28}) holds by Lemma (\ref{lem:2}) statement (2), Lemma (\ref{lem:5}) statement (2) and Lemma (\ref{lem:6}) statement (2). In the same way as Theorem 1.2 in \cite{Levendorski} we obtain defining relations (\ref{rel:30}), (\ref{rel:27}) and (\ref{rel:32}). Thus we finished the proof.

\begin{lemma} \label{lem:7}
\begin{enumerate}
\item The relations (\ref{rel:30}), (\ref{rel:27}) hold in $Y_{\hbar}^1(\hat{sl}(E,\Pi,p))$ when $\alpha_i\neq \alpha_j$.
\item The relation (\ref{rel:32}) holds for every roots $\alpha_i, \alpha_j$ in $Y_{\hbar}^1(\hat{sl}(E,\Pi,p))$.
\end{enumerate}
\end{lemma}

Relation (\ref{rel:27}) holds by Lemma (\ref{lem:5}) statement (1), Lemma (\ref{lem:6}) statement (5) and Lemma (\ref{lem:7}) statement (1). Relation (\ref{rel:30}) holds by Lemma (\ref{lem:5}) statement (4), Lemma (\ref{lem:6}) statement (4) and Lemma (\ref{lem:7}) statement (1).\\
We obtained (\ref{rel:33}), since (\ref{rel:33}) is equivalent to (\ref{rel:31}) when root is odd.Thus we need to show that relation (\ref{rel:34}) holds.

\begin{lemma}
\label{lem:8}
The relation (\ref{rel:34}) holds for odd roots in $Y_{\hbar}^1(\hat{sl}(E,\Pi,p))$.\\
\end{lemma}

We prove this relation in a similar way as done in \cite{Mazurenko}:
Let $X^{\pm}(r,0,0,s)$ be the left hand side of (\ref{rel:34}). We prove this relation by induction on $r$ and $s$ $\in \mathbb{Z}_{+}$. The initial case then $(r,0,0,s)=(0,0,0,0)$ is our initial assumption. If apply $ad(\tilde{h_{\alpha_{m},1}})$, $ad(\tilde{h_{\alpha_{n},1}})$ and $ad(\tilde{h_{\alpha_{k},1}})$ to $X^{\pm}(r,0,0,s)$ by relation (\ref{rel:40}) we obtain
\[
 a_{m,j-1}X^{\pm}(r+1,0,0,s)+a_{m,j}X^{\pm}(r,1,0,s)+X^{\pm}(r,0,1,s)+a_{m,j+1}X^{\pm}(r,0,0,s+1)=0,\\
\]
\[
 a_{n,j-1}X^{\pm}(r+1,0,0,s)+a_{n,j}X^{\pm}(r,1,0,s)+X^{\pm}(r,0,1,s)+a_{n,j+1}X^{\pm}(r,0,0,s+1)=0,\\
\]
\[
 a_{k,j-1}X^{\pm}(r+1,0,0,s)+a_{k,j}X^{\pm}(r,1,0,s)+X^{\pm}(r,0,1,s)+a_{k,j+1}X^{\pm}(r,0,0,s+1)=0.\\
\]
Since we have $A(m|n)$ type of algebra with $m,n\geq 2$ we obtain that $\alpha_{j-1}\neq \alpha_{j+1}$.
Consider the Cartan matrix block
\[
\hat{A}=
\begin{pmatrix}
  a_{m,j-1} & a_{m,j} & a_{m,j+1}\\
  a_{n,j-1} & a_{n,j} & a_{n,j+1}\\
  a_{k,j-1} & a_{k,j} & a_{k,j+1}\\
\end{pmatrix}.
\]
In order to determine when the determinant of $\hat{A}$ is nonzero consider the following the Dynkin diagrams(in other cases determinant of $\hat{A}$ will be zero):
\begin{enumerate}
    \item $\alpha_{j-1},\alpha_{j+1}$ are odd $\Rightarrow$ $m=j-2$, $n=j$, $k=j+1$,
    \item $\alpha_{j-1}$ is even $\alpha_{j+1}$ is odd $\Rightarrow$ $m=j-1$, $n=j$, $k=j+1$,
    \item $\alpha_{j-1}$ is odd $\alpha_{j+1}$ is even $\Rightarrow$ $m=j-1$, $n=j$, $k=j+1$,
    \item $\alpha_{j-1}$, $\alpha_{j+1}$ are even $\Rightarrow$ $m=j-2$, $n=j$, $k=j+1$.
\end{enumerate}
\begin{enumerate}
    \item $\alpha_{j-1},\alpha_{j+1}$ are odd $\Rightarrow$ $m=j+2$, $n=j$, $k=j+1$,
    \item $\alpha_{j-1}$ is even $\alpha_{j+1}$ is odd $\Rightarrow$ $m=j-1$, $n=j$, $k=j+1$,
    \item $\alpha_{j-1}$ is odd $\alpha_{j+1}$ is even $\Rightarrow$ $m=j-1$, $n=j$, $k=j+1$,
    \item $\alpha_{j-1}$, $\alpha_{j+1}$ are even $\Rightarrow$ $m=j+2$, $n=j$, $k=j+1$.
\end{enumerate}
When the determinant is nonzero we have $X^{\pm}(r+1,0,0,s)=X^{\pm}(r,0,0,s+1)=0$. The result follows by induction hypothesis.
Thus we proved the lemma.

Thus we have shown that (\ref{rel:34}) holds. We obtained relations (\ref{rel:33}) and (\ref{rel:34}) are deduced from relations (\ref{rel:35})-(\ref{rel:43}). This completes the proof of Theorem 1.\\
We also obtained minimalistic presentation of $Y_{\hbar}^1(\hat{sl}(E,\Pi,p))$

\begin{theorem}
 Suppose that $sl(E,\Pi,p)$ is the $A(m|n)$ type of algebra with $m\neq n$ and $m,n \geq 2$. Then $Y_{\hbar}^1(\tilde{sl}(E,\Pi,p))$ is isomorphic to the superalgebra generated by $x_{\alpha_i,r}^{\pm}, h_{\alpha_i,r}$ for $r \in \{0,1\}$ subject to defining relations (\ref{rel:35})-(\ref{rel:43}) and
\begin{equation}
\label{rel:79}
    [d,h_{\alpha_{\alpha_i},r}]=0,\quad  [d, x_{\alpha_i,r}^{+}]=
    \begin{cases}
    x_{\alpha_i,r}^{+}, & if \quad i=0\\
    0, & if \quad i \neq 0
    \end{cases}
    \quad
    [d,x_{\alpha_i,r}^{-}]=
    \begin{cases}
    - x_{\alpha_i,r}^{-}, & if \quad i=0,\\
    0, & if \quad i \neq 0.
    \end{cases}
\end{equation}
where odd generators $x_{\alpha_i,0}^{\pm}$ are corresponds odd roots and other generators are even.\\
\end{theorem}

The relation (\ref{in:1}) derived from relation (\ref{rel:79}) in a similar way to the Lemma (\ref{lem:1}) We proved the statement.

\subsection{Proof of Theorem \ref{thm_3.2}}

Since all elements of groupoid $\hat{W}$ can be represented as product of simple reflections $s_{\alpha_i}$. It is enough to check the isomorphism for every simple reflection $T_{s_{\alpha_i}}:Y(\hat{sl}(E,\Pi,p))\rightarrow Y(\hat{sl}(E,\Pi_1,p)$.
\begin{lemma}
\label{lem:9}
There exists isomorphism between of superalgebras $sl(E,\Pi,p)\rightarrow sl(E, \Pi_1, p)$ where map $\Pi\rightarrow \Pi_1$ induced by simple reflection $s_{\alpha_2}:\Pi\rightarrow \Pi_1$.\\
\end{lemma}
We denote root $\alpha_i$ as $\alpha_2$ to simplify the notation.
Since superalgebra $sl(E,\Pi,p)$ in distinguished realization is type $A(m|n)$. System of simple roots $\Delta$ does not change except roots $\alpha_1$,$\alpha_2$ and $\alpha_3$. By (\ref{rel:5})
\begin{enumerate}
    \item $s_{\alpha_2}:\alpha_1 \rightarrow \alpha_1+\alpha_2=\beta_1$,
    \item $s_{\alpha_2}:\alpha_2 \rightarrow -\alpha_2=\beta_2$,
    \item $s_{\alpha_2}:\alpha_3 \rightarrow \alpha_2+\alpha_3=\beta_3$.
\end{enumerate}

If the root $\alpha_2$ is even, the root system $\Pi$ of maps to itself $s_{\alpha_2}:\Pi \rightarrow \Pi$. Since the root $\alpha_2$ is even the parity of roots $\alpha_1,\alpha_3$ do not change, thus image of parity function on simple roots does not change. The map of root systems induces map of super algebras $\hat{s}_{\alpha_2}:sl(E,\Pi,p)\rightarrow sl(E, \Pi, p)$. We introduce it as in Kodera work \cite{Kodera}:
\[
\hat{s}_{\alpha_i}(x_{\alpha_j}^{\pm})=
\begin{cases}
         -x_{\alpha_i}^{\mp} & if \quad a_{i,j}=2,\\
        \pm[x_{\alpha_i}^{\pm},x_{\alpha_j}^{\pm}], & if \quad a_{i,j}=-1,\\
         x_{\alpha_j}^{\pm} & if \quad a_{i,j}=0,\\
\end{cases}
\quad
\hat{s}_{\alpha_i}(h_{\alpha_j,0})=
\begin{cases}
    -h_{\alpha_i}, & if \quad a_{i,j}=2,\\
    h_{\alpha_i}+h_{\alpha_j}, & if \quad a_{i,j}=-1, \\
    h_{\alpha_j}, & if \quad a_{i,j}=0.
\end{cases}
\]

We introduce map of Yangians $T_{s_{\alpha_2}}$ that coincides on the generators of zero order with the map $\hat{s}_{\alpha_2}$.
We define map of Yangians $s_{\alpha_2}$ for generators of non zero order as in Kodera work \cite{Kodera}:
\[
T_{s_{\alpha_i}}(x_{\alpha_j,1}^{\pm})=
\begin{cases}
    -x_{\alpha_i,1}^{\mp}+\frac{\hbar}{2}\{h_{\alpha_i,0},x_{\alpha_i,0}^{\mp}\}, & if \quad a_{i,j}=2,\\
    \pm[x_{\alpha_i,0}^{\pm},x_{\alpha_j,1}^{\pm}], & if \quad a_{i,j}=-1,\\
     x_{\alpha_j,1}^{\pm}, & if \quad a_{i,j}=0.\\
\end{cases} \quad
\]

\[ T_{s_{\alpha_i}}(\tilde{h}_{\alpha_j,1})=
\begin{cases}
    -\tilde{h}_{\alpha_j,1}-\hbar\{x_{\alpha_i,0}^{+},x_{\alpha_i,0}^{-}\},& if \quad a_{i,j}=2,\\
    \tilde{h}_{\alpha_j,1}+\tilde{h}_{\alpha_i,1}+\frac{\hbar}{2}\{x_{\alpha_i,0}^{+},x_{\alpha_i,0}^{-}\}, & if \quad a_{i,j}=-1, \\
    \tilde{h}_{\alpha_j,1},& if \quad a_{i,j}=0.\\
\end{cases}
\]

If the root $\alpha_2$ is odd, the image of parity function of $\alpha_1, \alpha_3$ changes: $p(\beta_1)=p(\alpha_1)+1$, $p(\beta_3)=p(\alpha_3)+1$, since the image of function $p(\Pi)$ lies in $\mathbb{Z}_2$. Thus the corresponding block of Cartan matrix changes. Since the roots $\alpha_1, \alpha_3$ changing in similar way, it is enough to check the defining relations for $\alpha_3$ and $\alpha_2$.
The block is changing the following way, we prove the case when $\alpha_1=\delta_1-\delta_2$ since the proof when $\alpha_1=\varepsilon_1-\varepsilon_2$ is the same. By definition of Cartan matrix we have $(\alpha_1,\alpha_1)=2$ and $(\alpha_1,\alpha_2)=-1$, hence $(-\alpha_2,-\alpha_2)=0=(\beta_2,\beta_2)$,$(\alpha_1+\alpha_2,-\alpha_2)=-1=(\beta_1,\beta_2)$ and $(\alpha_1+\alpha_2,\alpha_1+\alpha_2)=0=(\beta_1,\beta_1)$. Thus, we obtain the structure of block of Cartan matrix for $sl(\Pi_1)$
\begin{equation}
\begin{pmatrix}
     2 & -1\\
    -1 & 0
    \end{pmatrix} \rightarrow \begin{pmatrix}
    0 & 1\\
    1 & 0
    \end{pmatrix}.
\end{equation}

We introduce the map $\phi:sl(\Pi)\rightarrow sl(\Pi_1)$, which induced by an element of Weyl groupoid.  We show that this map is Lie superalgebra isomorphism. This map is composition of the maps which induced by reflection. It suffices to prove that this mapping is an isomorphism in the case when it is induced by one reflection. It is well known that this map is isomorphism, when it is induced by even reflection.  Then it is sufficient to prove this assertion in the case, when  $\phi:sl(\Pi)\rightarrow sl(\Pi_1)$ is induced by odd reflection $s_i$. In this case $\psi$ it is defined on root generators  by the following formulas:


\begin{enumerate}
    \item $x_{\beta_1}^{\pm}=\pm[x_{\alpha_1}^{\pm},x_{\alpha_2}^{\pm}]$,
    \item $x_{\beta_2}^+=x_{\alpha_2}^-$,
    \item $x_{\beta_2}^-=x_{\alpha_2}^+$,
    \item $x_{\beta_3}^{\pm}=\pm[x_{\alpha_3}^{\pm},x_{\alpha_2}^{\pm}]$,
    \item $h_{\beta_1}=h_{\alpha_1}+h_{\alpha_2}$,
    \item $h_{\beta_3}=h_{\alpha_3} + h_{\alpha_2}$,
    \item $h_{\beta_2}=-h_{\alpha_2}$.
\end{enumerate}

We can define quantum reflections with respect to the odd simple root by the following formula. 

\[
T_{s_{\alpha_i}}(x_{\alpha_j,1}^{\pm})=
\begin{cases}
        \pm[x_{\beta_i,0}^{\pm},x_{\beta_j,1}^{\pm}], & if \quad a_{i,j}=-1,\\
        x_{-\beta_i,1}^{\pm}, & if \quad i = j, a_{i,j}=0.\\
     x_{\alpha_j,1}^{\pm}, & if \quad i\neq j, a_{i,j}=0.\\
\end{cases} \quad
\]

\[ T_{s_{\alpha_i}}(\tilde{h}_{\alpha_j,1})=
\begin{cases}
        \tilde{h}_{\beta_j,1}+\tilde{h}_{\beta_i,1}+\frac{\hbar}{2}\{x_{\beta_i,0}^{+},x_{\beta_i,0}^{-}\}, & if \quad a_{i,j}=-1, \\
 \tilde{h}_{-\beta_j,1},& if \quad i=j.\\
    \tilde{h}_{\alpha_j,1},& if \quad i\neq j, a_{i,j}=0.\\
\end{cases}
\]

We want to prove the relation for realization $sl(\Pi_1)$. We prove relations (\ref{in:6}) only for block of Cartan matrix corresponding to roots $\alpha_1,\alpha_2$, since $h_{\alpha_3}$ changes similar to $\alpha_{1}$. First relation is obvious since by structure of map $\phi$ and first relation (\ref{in:6}) for realization $sl(\Pi)$: $[h_{\beta_1},h_{\beta_2}]\rightarrow [h_{\alpha_1}+h_{\alpha_2},-h_{\alpha_2}]=0$. The relations $[h_{\beta_1},h_{\beta_1}]=0$ and $[h_{\beta_2},h_{\beta_2}]=0$ proved in similar way.\\
By second relation from (\ref{in:6}) we prove relation $[h_{\beta_1},x_{\beta_1}^+]\rightarrow[h_{\alpha_1}+h_{\beta_2},[x_{\alpha_1}^+,x_{\alpha_2}^+]] = 2[x_{\alpha_1}^+,x_{\alpha_2}^+]-[x_{\alpha_1}^+,x_{\alpha_2}^+] - [x_{\alpha_1}^+,x_{\alpha_2}^+]=0$, relations $[h_{\beta_1},x_{\beta_2}^+]$,$[h_{\beta_2},x_{\beta_2}^+]$ and $[h_{\beta_1},x_{\beta_2}^+]$ proves in a similar way. Thus the second relation from (\ref{in:6}) is proved.
By the structure of map $\phi$ and third relation (\ref{in:6}) we prove the following relation
\begin{equation}\label{4.36}
[x_{\beta_1}^-,x_{\beta_2}^+]\rightarrow [[[x_{\alpha_1}^-,x_{\alpha_2}^-],[x_{\alpha_1}^-,x_{\alpha_2}^-]]= [[x_{\alpha_1}^-,x_{\alpha_2}^-],x_{\alpha_1}^+],x_{\alpha_2}^+] + [x_{\alpha_1}^-,[[x_{\alpha_1}^-,x_{\alpha_2}^-],x_{\alpha_2}^+].
\end{equation}
By third relation (\ref{in:6}) we obtain that (\ref{4.36}) is equal to
\begin{equation}
    -[[h_{\alpha_1},x_{\alpha_2}^-], x_{\alpha_2}^+]-[x_{\alpha_1}^+, [x_{\alpha_1}^-,h_{\alpha_2}]].
\end{equation}
Using second and third relation (\ref{in:6}) for realization $sl(\Pi)$ we obtain
\begin{equation}
    -h_{\alpha_2}-h_{\alpha_1}.
\end{equation}

The relation $[x_{\beta_2}^+,x_{\beta_2}^-] = \delta_{ij}h_{\alpha_2}$ directly follows from structure of map $\phi$ and third relation (\ref{in:6}) for realization $sl(\Pi)$. By Serre relation $ad(x_{j})^{1+|a_{i,j}|}$ for realization $sl(\Pi)$, since $|a_{12}|=1$ we prove relations $[x_{\beta_1}^+,x_{\beta_2}^-]=0$ and $[x_{\beta_2}^+,x_{\beta_1}^-]=0$, hence we proved third relation (\ref{in:6}) for $sl(\Pi_1)$.\\
The relation $ad(x_{\beta_1}^{1+|a_{ij}|}(x_{\beta_2}))=0$ and relation $ad(x_{\beta_1}^{1+|a_{ij}|}(x_{\beta_2}))=0$ follows directly from Serre relations, hence we proved Serre relations for $sl(\Pi_1)$\\
By structure of map $\phi$ we obtain
\begin{equation}
\label{ref:91}
    [x_{\beta_1}^{\pm},x_{\beta_1}^{\pm}]\rightarrow [[x_{\alpha_1}^{\pm},x_{\alpha_2}^{\pm}],[x_{\alpha_1}^{\pm},x_{\alpha_2}^{\pm}]],
\end{equation}
\begin{equation}
\label{rel:92}
    [x_{\beta_2}^{\pm},x_{\beta_2}^{\pm}]\rightarrow [x_{\alpha_2}^{\mp},x_{\alpha_2}^{\mp}]=0.
\end{equation}
Relation \ref{rel:92} holds by first relation (\ref{in:7}).
By Serre relations for $sl(\Pi)$ we obtain following equation from (\ref{ref:91})
\begin{equation}
[[[x_{\alpha_1}^+,x_{\alpha_2}^+],x_{\alpha_1}^+],x_{\alpha_2}^+] + [x_{\alpha_1}^+,[[x_{\alpha_1}^+,x_{\alpha_2}^+],x_{\alpha_2}^+]] = 0.
\end{equation}
Thus we obtained first relation from (\ref{in:7})
By the structure of map $\phi$ we obtain
\begin{equation}
\label{rel:94}
[[x_{\beta_1}^{\pm},x_{\beta_2}^{\pm}],[x_{\beta_2}^{\pm},x_{\beta_3}^{\pm}]]\rightarrow [[[x_{\alpha_1}^{\pm},x_{\alpha_2}^{\pm}],x_{\alpha_2}^{\mp}],[,x_{\alpha_2}^{\mp},[x_{\alpha_2}^{\pm},x_{\alpha_3}^{\pm}]]].
\end{equation}
By Serre relation for $sl(\Pi)$ we obtain that first part of right hand side of following relation for image of map $\phi$ in (\ref{rel:94}) is equal zero
\begin{equation}
\label{rel:95}
[[[x_{\alpha_1}^{\pm},x_{\alpha_2}^{\pm}],x_{\alpha_2}^{\mp}],[,x_{\alpha_2}^{\mp},[x_{\alpha_2}^{\pm},x_{\alpha_3}^{\pm}]]]= [[[x_{\alpha_1}^{\pm},x_{\alpha_2}^{\pm}],x_{\alpha_2}^{\mp}],x_{\alpha_2}^{\mp}]+ [x_{\alpha_2}^{\mp},[[x_{\alpha_1}^{\pm},x_{\alpha_2}^{\pm}],[x_{\alpha_2}^{\pm},x_{\alpha_3}^{\pm}]]] = 0.
\end{equation}
The second part of right hand side of (\ref{rel:95}) is equal zero by second relation (\ref{in:7}) in realization $sl(\Pi)$.
By the structure of map $\phi$ we obtain
\begin{equation}
\label{rel:96}
[x_{\beta_0}^{\pm},x_{\beta_1}^{\pm},[x_{\beta_1}^{\pm},x_{\beta_2}^{\pm}]]\rightarrow [[[x_{\alpha_0}^{\pm},[x_{\alpha_1}^{\pm},x_{\alpha_2}^{\pm}]],[[x_{\alpha_1}^{\pm},x_{\alpha_2}^{\pm}],x_{\alpha_2}^{\mp}]]].
\end{equation}
Where the root $\alpha_0$ is the root next to root $\alpha_1$ on Dynkin diagram, by the definition of simple reflection $s_{\alpha_2}(\alpha_0)=\alpha_0$, hence $\phi(x_{\alpha_0}^{\pm})=x_{\beta_0}^{\pm}$.
By third relation in (\ref{in:6}) we obtain
\begin{equation}
\label{rel:97}
[[[x_{\alpha_0}^{\pm},[x_{\alpha_1}^{\pm},x_{\alpha_2}^{\pm}]],[[x_{\alpha_1}^{\pm},x_{\alpha_2}^{\pm}],x_{\alpha_2}^{\mp}]]]= [[x_{\alpha_0}^{\pm},[x_{\alpha_1}^{\pm},x_{\alpha_2}^{\pm}]],[x_{\alpha_1}^{\pm},h_{\alpha_2}]]=
-[[x_{\alpha_0}^{\pm},[x_{\alpha_1}^{\pm},x_{\alpha_2}^{\pm}]],x_{\alpha_1}^{\pm}].
\end{equation}
Where the second equation follows from (\ref{in:6}) second relation for $sl(\Pi)$
By second relation in (\ref{in:7}) and Serre relations for $sl(\Pi)$ we obtain following relation
\begin{equation}
\label{rel:98}
-[[x_{\alpha_0}^{\pm},[x_{\alpha_1}^{\pm},x_{\alpha_2}^{\pm}]],x_{\alpha_1}^{\pm}]= -(2[[x_{\alpha_0}^{\pm},x_{\alpha_1}^{\pm}],[x_{\alpha_1}^{\pm},x_{\alpha_2}^{\pm}]]+ [[x_{\alpha_0}^{\pm},x_{\alpha_1}^{\pm}],[x_{\alpha_1}^{\pm}],x_{\alpha_2}^{\pm}])+ [[x_{\alpha_1}^{\pm},[x_{\alpha_1}^{\pm},x_{\alpha_0}^{\pm}]],x_{\alpha_2}^{\pm}] = 0.
\end{equation}
For roots $\alpha_2,\alpha_3,\alpha_4$ proof is similar to proof for roots $\alpha_0,\alpha_1,\alpha_2$ by relations (\ref{rel:96}), (\ref{rel:97}), (\ref{rel:98}). Thus we proved all relations (\ref{in:6}), (\ref{in:7}) for superalgebra $sl(\Pi_1)$.\\
The inverse map $\phi^{-1}$ is induced by simple reflection $s_{\beta_2}$. Note that map $\phi^{-1}$ translates system of simple roots $\Pi_1$ to $\Pi$, hence corresponding block of Cartan matrix changes as in following map
\begin{equation}
   \begin{pmatrix}
     0 & 1\\
    1 & 0
    \end{pmatrix} \rightarrow \begin{pmatrix}
    2 & -1\\
    -1 & 0
    \end{pmatrix}
\end{equation}
The defining relations proves similar to relations for map $\phi$.
Since any element of Weyl groupoid $\hat{W}$ is equal to product of elements corresponding to simple reflections, we can construct sequence of isomorphisms of $sl(\Pi_i)$ superalgebras for action of any element of Weyl groupoid, hence construct isomorphism for action of any element of Weyl groupoid $\hat{W}$.
We obtained the proof of Lemma (\ref{lem:9})\\

\subsection{Proof of theorem 3.3}

We construct isomorphism $T_s:Y(\hat{sl}(E,\Pi,p))\rightarrow Y(\hat{sl}(E,\Pi_1,p)$
as in Lemma (\ref{lem:2}) we present element of Weyl groupoid $s \in \hat{W}$ as product of elements $s_{\alpha_i}$ corresponding to simple reflections $s_{\alpha_i}$.\\
We rename our roots as in Lemma (\ref{lem:2}) to simplify the notation.\\
In Theorem 3.1 we proved that Yangian $Y(sl(\Pi))$ of arbitrary realization is isomorphic to Yangian in minimalistic system of generators. We will use this fact to construct a map.
We define map $\psi:Y(sl(\Pi))\rightarrow Y(sl(\Pi_1))$, which induced by element of Weyl groupoid. More precisely, $\psi = T_w := T_{s_1}\circ \ldots \circ T_{s_k}$, where $w = s_1 \ldots s_k \in W$ is an element of Weyl groupoid.  Easy to see that general case can be reduced to the partial case when $\psi= T_{s_i}$.  We also will consider root generators of $Y_{\hbar}(\hat{sl}(m|n))$ as images of root generators of $Y_{\hbar}(sl(1|2))$ under natural embedding $Y_{\hbar}(sl(1|2)) \rightarrow Y_{\hbar}(\hat{sl}(m|n))$.  In this case when $\psi= T_{s_i}$ is induced by odd reflection $s_i$ it is defined on root generators  by the following formulas:
\begin{enumerate}
\item $x_{\beta_1,0}^{\pm}=\pm[x_{\alpha_1,0}^{\pm},x_{\alpha_2,0}^{\pm}]$,
\item $x_{\beta_2,0}^+=x_{\alpha_2,0}^-$,
\item $x_{\beta_2,0}^-=x_{\alpha_2,0}^+$,
\item $x_{\beta_3,0}^{\pm}=\pm[x_{\alpha_3,0}^{\pm},x_{\alpha_2,0}^{\pm}]$,
\item $\tilde{h}_{\beta_1}=\tilde{h}_{\alpha_1}+\tilde{h}_{\alpha_2}+ \frac{\hbar}{2}\{x_{\alpha_{2},0}^{+},x_{\alpha_{2},0}^{-}\}$,
\item $\tilde{h}_{\beta_3}=\tilde{h}_{\alpha_3}+\tilde{h}_{\alpha_2}+ \frac{\hbar}{2}\{x_{\alpha_{2},0}^{+},x_{\alpha_{2},0}^{-}\}$,
\item $\tilde{h}_{\beta_2}=-\tilde{h}_{\alpha_2}$.
\end{enumerate}


Note that on the elements of zero-grading Yangian $Y(sl(\Pi))$ has structure of $sl(\Pi)$, hence relations
(\ref{rel:35}), (\ref{rel:41}), (\ref{rel:42}), (\ref{rel:43}) is proved by Lemma 4.1.2. Thus we only need to prove relations (\ref{rel:37}), (\ref{rel:39}), (\ref{rel:40}).
\begin{lemma}
\label{lemma:4.10}
    The relation $(\ref{rel:39})$ holds for map $\psi$.
\end{lemma}

By definition of map $\psi$ and relation $\ref{rel:40}$ we can obtain
\begin{equation}
[x_{\beta_1,1}^+,x_{\beta_2,0}^+] = -[[\tilde{h}_{\beta_2,1},x_{\beta_1,0}^+],x_{\beta_2,0}^+] \rightarrow [[[\tilde{h}_{\alpha_2,1},[x_{\alpha_1,0}^+,x_{\alpha_2,0}^+]],x_{\alpha_2,0}^-].
\end{equation}
By relation (\ref{rel:41}) we obtain
\begin{equation}
\label{rel:100}
[[x_{\alpha_1,1}^+,x_{\alpha_2,0}^+],x_{\alpha_2,0}^-]= [[x_{\alpha_1,1}^+,x_{\alpha_2,0}^-],x_{\alpha_2,0}^+]+[x_{\alpha_1,1}^+,[x_{\alpha_2,0}^+,x_{\alpha_2,0}^-]].
\end{equation}
First term of right hand side of (\ref{rel:100}) vanishes by the relation (\ref{rel:37}), by the relation (\ref{rel:36}) we obtain the following relation from (\ref{rel:100})
\begin{equation}
\label{rel:98}
    [x_{\alpha_1,1}^+, h_{\alpha_2,0}] = -x_{\alpha_1,1}^+.
\end{equation}
Thus we obtained
\begin{equation}
\label{midres:1}
[x_{\beta_1,1}^+,x_{\beta_2,0}^+]\rightarrow -x_{\alpha_1,1}^+.
\end{equation}
By definition of map $\psi$ and relation (\ref{rel:40}) we can obtain
\begin{equation}
\label{equ:100}
[x_{\beta_1,0}^+,x_{\beta_2,1}^+]=-[x_{\beta_1,0}^+,[\tilde{h}_{\beta_1,1},x_{\beta_2,0}]]\rightarrow -[[x_{\alpha_1,0}^+,x_{\alpha_2,0}^+],[\tilde{h}_{\alpha_1,1}+\tilde{h}_{\alpha_2,1}+ \frac{\hbar}{2}\{x_{\alpha_{2},0}^{+},x_{\alpha_{2},0}^{-}\},x_{\alpha_2,0}^-]].
\end{equation}
By applying relation (\ref{rel:40}) to the right hand side of (\ref{equ:100}) we obtain
\begin{equation}
\label{equ:101}
\begin{split}
&-[[x_{\alpha_1,0}^+,x_{\alpha_2,0}^+],[\tilde{h}_{\alpha_1,1}+\tilde{h}_{\alpha_2,1}]
 +\frac{\hbar}{2}\{x_{\alpha_{2},0}^{+},x_{\alpha_{2},0}^{-}\},x_{\alpha_2,0}^-]]=\\
&=[[x_{\alpha_{1},0}^{+},x_{\alpha_{2},0}^{+}],x_{\alpha_{2},1}^{-}]- \frac{\hbar}{2}[[x_{\alpha_{1},0}^{+},x_{\alpha_{2},0}^{+}],[\{x_{\alpha_{2},0}^{+},x_{\alpha_{2},0}^{-},x_{\alpha_{2},0}^{-}\}]].
\end{split}
\end{equation}
By relation (\ref{rel:40}) we obtain
\begin{equation}
[[x_{\alpha_1,0}^+,x_{\alpha_2,0}^+],x_{\alpha_2,1}^-]= [[x_{\alpha_1,0}^+,x_{\alpha_2,0}^-],x_{\alpha_2,1}^+]+[x_{\alpha_1,0}^+,[x_{\alpha_2,0}^+,x_{\alpha_2,0}^-]],
\end{equation}
where first term of right hand side  is equal to zero by (\ref{rel:36}), second term of right hand side by (\ref{rel:36}) equals to $[x_{\alpha_1,0}^+,h_{\alpha_2,1}]$. By definition of $\tilde{h}_{\alpha_i,1}$
\begin{equation}
\label{rel:101}
[x_{\alpha_i,0}^+,h_{\alpha_2,1}]=[x_{\alpha_i,0}^+,\tilde{h}_{\alpha_2,1}+ \frac{\hbar}{2}h_{\alpha_2,0}^2]=-x_{\alpha_1,1}^++\frac{\hbar}{2}\{h_{\alpha_2,0},x_{\alpha_2,0}^+\}.
\end{equation}

By definition of map $\psi$
\begin{equation}
\label{midres:3}
\frac{\hbar}{2}\left(x_{\beta_{1},0}^{+}x_{\beta_{2},0}^{+}+x_{\beta_{2},0}^{+}x_{\beta_{1},0}^{+}\right) \rightarrow \frac{\hbar}{2}\{[x_{\alpha_{1},0}^{+};x_{\alpha_{2},0}^{+}];x_{\alpha_{2},0}^{-}\}.
\end{equation}
By the Jacobi identity we have following relation for second term of right hand side of (\ref{equ:101})
\begin{equation}
\label{equ:105}
\begin{split}
&\frac{\hbar}{2}[[x_{\alpha_{1},0}^{+},x_{\alpha_{2},0}^{+}],[x_{\alpha_{2},0}^{+}x_{\alpha_{2},0}^{-}+ x_{\alpha_{2},0}^{-}x_{\alpha_{2},0}^{+},x_{\alpha_{2},0}^{-}]]= [[x_{\alpha_{1},0}^{+},x_{\alpha_{2},0}^{+}],[x_{\alpha_{2},0}^{+},x_{\alpha_{2},0}^{-}]x_{\alpha_{2},0}^{-}+\\
&+x_{\alpha_{2},0}^{+}[x_{\alpha_{2},0}^{-},x_{\alpha_{2},0}^{-}]
+[x_{\alpha_{2},0}^{-},x_{\alpha_{2},0}^{-}]x_{\alpha_{2},0}^{+} + x_{\alpha_{2},0}^{-}[x_{\alpha_{2},0}^{+},x_{\alpha_{2},0}^{-}]].
\end{split}
\end{equation}
By applying relation (\ref{rel:36}) to the right hand side of (\ref{equ:105}) we obtain
\begin{equation}
\label{equ:106}
\begin{split}
&\frac{\hbar}{2}[[x_{\alpha_{1},0}^{+},x_{\alpha_{2},0}^{+}],[x_{\alpha_{2},0}^{+},x_{\alpha_{2},0}^{-}]x_{\alpha_{2},0}^{-}+ x_{\alpha_{2},0}^{+}[x_{\alpha_{2},0}^{-},x_{\alpha_{2},0}^{-}]
+ [x_{\alpha_{2},0}^{-},x_{\alpha_{2},0}^{-}]x_{\alpha_{2},0}^{+}+ x_{\alpha_{2},0}^{-}[x_{\alpha_{2},0}^{+},x_{\alpha_{2},0}^{-}]]=\\
& = \frac{\hbar}{2}[[x_{\alpha_{1},0}^{+},x_{\alpha_{2},0}^{+}],h_{\alpha_{2},0}x_{\alpha_2,0}^{-}+ 0+ 0+ x_{\alpha_2,0}^{+}h_{\alpha_{2},0}].
\end{split}
\end{equation}
By the Jacobi identity from equation (\ref{equ:106}) we obtain
\begin{equation}
\label{equ:107}
\begin{split}
&\frac{\hbar}{2}[[x_{\alpha_{1},0}^{+},x_{\alpha_{2},0}^{+}],h_{\alpha_{2},0}x_{\alpha_2,0}^{-}+ x_{\alpha_2,0}^{+}h_{\alpha_{2},0}] = \frac{\hbar}{2}([[x_{\alpha_{1},0}^{+},x_{\alpha_{2},0}^{+}],h_{\alpha_{2},0}]x_{\alpha_{2},0}^{-} + h_{\alpha_{2},0}[[x_{\alpha_{1},0}^{+},x_{\alpha_{2},0}^{+}],x_{\alpha_{2},0}^{-}]+\\
& + [[x_{\alpha_{1},0}^{+},x_{\alpha_{2},0}^{+}],x_{\alpha_{2},0}^{-}]h_{\alpha_{2},0} + x_{\alpha_{2},0}^{-}[[x_{\alpha_{1},0}^{+},x_{\alpha_{2},0}^{+}],h_{\alpha_{2},0}]).
\end{split}
\end{equation}
By Jacobi identity and applying relations (\ref{rel:36}), (\ref{rel:42}) and (\ref{rel:38}) to the right hand side of we obtain following equation
\begin{equation}
\begin{split}
&\frac{\hbar}{2}([[x_{\alpha_{1},0}^{+}, x_{\alpha_{2},0}^{+}], h_{\alpha_{2},0}]x_{\alpha_{2},0}^{-}+h_{\alpha_{2},0}[[x_{\alpha_{1},0}^{+},x_{\alpha_{2},0}^{+}],x_{\alpha_{2},0}^{-}] +
[[x_{\alpha_{1},0}^{+},x_{\alpha_{2},0}^{+}],x_{\alpha_{2},0}^{-}]h_{\alpha_{2},0} +\\
&+x_{\alpha_{2},0}^{-}[[x_{\alpha_{1},0}^{+},x_{\alpha_{2},0}^{+}],h_{\alpha_{2},0}]) =
\frac{\hbar}{2}(-[x_{\alpha_{1},0}^{+},x_{\alpha_{2},0}^{+}]x_{\alpha_{2},0}^{-} - x_{\alpha_{2},0}^{-}[x_{\alpha_{1},0}^{+},x_{\alpha_{2},0}^{+}] + h_{\alpha_{2},0}x_{\alpha_{1},0}+x_{\alpha_{1},0}h_{\alpha_{2},0}) = \\
&= \frac{\hbar}{2}(-\{[x_{\alpha_{1},0}^{+};x_{\alpha_{2},0}^{+}];x_{\alpha_{2},0}^{-}\} + \{h_{\alpha_{2},0}, x_{\alpha_{\alpha_{1},0}}\}).
\end{split}
\end{equation}

Combining this relation with (\ref{rel:101}),(\ref{equ:100}) and (\ref{equ:101}) we obtain
\begin{equation}
\label{midres:2}
\begin{split}
&[x_{\beta_1,0}^+,x_{\beta_2,1}^+]\rightarrow -x_{\alpha_1,1}^+ + \frac{\hbar}{2}\{h_{\alpha_2,0},x_{\alpha_2,0}^+\}+\frac{\hbar}{2}(\{[x_{\alpha_{1},0}^{+}; x_{\alpha_{2},0}^{+}];x_{\alpha_{2},0}^{-}\}) - \frac{\hbar}{2}\{h_{\alpha_{2},0},x_{\alpha_{\alpha_{1},0}}\} =\\
&= -x_{\alpha_1,1}^++\frac{\hbar}{2}\{[x_{\alpha_{1},0}^{+};x_{\alpha_{2},0}^{+}];x_{\alpha_{2},0}^{-}\}.
\end{split}
\end{equation}
Using (\ref{midres:1}) and (\ref{midres:2}) we obtain
\begin{equation}
\label{midresu:1}
[x_{\beta_1,1}^+,x_{\beta_2,0}^+]-[x_{\beta_1,0}^+,x_{\beta_2,1}^+]\rightarrow-\frac{\hbar}{2}\{[x_{\alpha_{1},0}^{+};x_{\alpha_{2},0}^{+}];x_{\alpha_{2},0}^{-}\}
\end{equation}
Using (\ref{midres:3}) and (\ref{midresu:1}) we prove that $\psi$ preserves relation (\ref{rel:39}).
Thus we proved the lemma \ref{lemma:4.10},

\begin{lemma}
\label{lemma:4.11}
    The relation $\ref{rel:37}$ holds for map $\psi$
\end{lemma}

By the definition of map $\psi$ and relation (\ref{rel:40}) we obtain
\begin{equation}
    [x_{\beta_{1},1}^{+},x_{\beta_{1},0}^{-}] = -[[\tilde{h}_{\beta_2,1},x_{\beta_1,0}^+], x_{\beta_1,0}^-]\rightarrow[[\tilde{h}_{\alpha_{2},1}, [x_{\alpha_{1},0}^{+},x_{\alpha_{2},0}^{+}]], [x_{\alpha_{1},0}^{-},x_{\alpha_{2},0}^{-}]].
\end{equation}
By the relation (\ref{rel:40}) we obtain
\begin{equation}
[[\tilde{h}_{\alpha_{2},1},[x_{\alpha_{1},0}^{+},x_{\alpha_{2},0}^{+}]],[x_{\alpha_{1},0}^{-},x_{\alpha_{2},0}^{-}]] = [[x_{\alpha_{1},1}^{+},x_{\alpha_{2},0}^{+}],[x_{\alpha_{1},0}^{-},x_{\alpha_{2},0}^{-}]].
\end{equation}
By the Jacobi identity we obtain
\begin{equation}
\label{equ:114}
\begin{split}
&[[x_{\alpha_{1},1}^{+},x_{\alpha_{2},0}^{+}],[x_{\alpha_{1},0}^{-},x_{\alpha_{2},0}^{-}]]=[[[x_{\alpha_{1},1}^{+},x_{\alpha_{2},0}^{+}],x_{\alpha_{1},0}^{-}],x_{\alpha_{2},0}^{-}]+[x_{\alpha_{1},0}^{-},[[x_{\alpha_{1},1}^{+},x_{\alpha_{2},0}^{+}]],x_{\alpha_{2},0}^{-}]]=\\
&=[[[x_{\alpha_{1},0}^{-},x_{\alpha_{2},0}^{+}],x_{\alpha_{1},1}^{+}]-[x_{\alpha_{2},0}^{+},[x_{\alpha_{1},0}^{-},x_{\alpha_{1},1}^{+}]],x_{\alpha_{2},0}^{-}]+[x_{\alpha_{1},0}^{-},[[x_{\alpha_{2},0}^{-},x_{\alpha_{2},0}^{+}],x_{\alpha_{1},1}]-\\
&-[x_{\alpha_{2},0}^{+},[x_{\alpha_{2},0}^{-}.x_{\alpha_{1},1}^{+}]]].
\end{split}
\end{equation}
By applying relations (\ref{rel:36}) and (\ref{rel:37}) to the right hand side of (\ref{equ:114}) we obtain
\begin{equation}
\label{equ:115}
\begin{split}
&[[[x_{\alpha_{1},0}^{-},x_{\alpha_{2},0}^{+}],x_{\alpha_{1},1}^{+}] - [x_{\alpha_{2},0}^{+}, [x_{\alpha_{1},0}^{-},x_{\alpha_{1},1}^{+}]], x_{\alpha_{2},0}^{-}] + [x_{\alpha_{1},0}^{-}, [[x_{\alpha_{2},0}^{-}, x_{\alpha_{2},0}^{+}], x_{\alpha_{1},1}^{+}] - \\
&-[x_{\alpha_{2},0}^{+},[x_{\alpha_{2},0}^{-}, x_{\alpha_{1},1}^{+}]]] = [0-[x_{\alpha_{2},0}^{+},-h_{\alpha_{1},1}], x_{\alpha_{2},0}^{-}]+[x_{\alpha_{1},0}^{-}, [-h_{\alpha_{2},0},x_{\alpha_{1},1}^{+}]-0].
\end{split}
\end{equation}
By applying relations (\ref{rel:38}) and (\ref{rel:37}) to the second term of right hand side of (\ref{equ:115})
we obtain
\begin{equation}
\label{equ:116}
[x_{\alpha_{1},0}^{-},[-h_{\alpha_{2},0},x_{\alpha_{2},1}^{+}]] = -[x_{\alpha_{1},0}^{-},x_{\alpha_{1},1}^{+}] = h_{\alpha_{1},1}.
\end{equation}
By the definition of $\tilde{h}_{\alpha_{1},1}$ and the relations (\ref{rel:40}),(\ref{rel:39}) and (\ref{rel:38}) to the first term of the right hand side of the (\ref{equ:115})  we obtain
\begin{equation}
\label{equ:117}
\begin{split}
&[[x_{\alpha_{2},0}^{+},-h_{\alpha_{1},1}],x_{\alpha_{2},0}^{-}] = -[[x_{\alpha_{2},0}^{+},\tilde{h}_{\alpha_{1},1}+\frac{\hbar}{2}h_{\alpha_{1},0}^{2}],x_{\alpha_{2},0}^{-}] = [x_{\alpha_{2},1}^{+},x_{\alpha_{2},0}^{-}] + \\
&+\frac{\hbar}{2}[\{h_{\alpha_{1},0},x_{\alpha_{2},0}^{+}\},x_{\alpha_{2},0}^{-}]=h_{\alpha_{2},1} + \frac{\hbar}{2}[\{h_{\alpha_{1},0},x_{\alpha_{2},0}^{+}\},x_{\alpha_{2},0}^{-}].
\end{split}
\end{equation}
By applying (\ref{rel:36}), (\ref{rel:38}) and Jacobi identity to the second term of the right hand side of (\ref{equ:117})
we obtain
\begin{equation}
\label{equ:118}
\begin{split}
&\frac{\hbar}{2}[\{h_{\alpha_{1},0},x_{\alpha_{2},0}^{+}\},x_{\alpha_{2},0}^{-}]=\frac{\hbar}{2}[h_{\alpha_{1},0}x_{\alpha_{2},0}^{+}+x_{\alpha_{2},0}^{+}h_{\alpha_{1},0},x_{\alpha_{2},0}^{-}]=\frac{\hbar}{2}([h_{\alpha_{1},0},x_{\alpha_{2},0}^{-}]x_{\alpha_{2},0}^{+}
+\\
&+h_{\alpha_{1},0}[x_{\alpha_{2},0}^{+},x_{\alpha_{2},0}^{-}]
+[x_{\alpha_{2},0}^{+},x_{\alpha_{2},0}^{-}]h_{\alpha_{1},0}+x_{\alpha_{2},0}^{+}[h_{\alpha_{1},0},x_{\alpha_{2},0}^{-}])= \frac{\hbar}{2}\{x_{\alpha_{2},0}^{+},x_{\alpha_{2},0}^{-}\}+\frac{\hbar}{2}\{h_{\alpha_{1},0},h_{\alpha_{2},0}\}.
\end{split}
\end{equation}
By the definition of $h_{\beta_{1},1}$, $h_{\alpha_{1},1}$ and map $\psi$ we obtain
\begin{equation}
\label{equ:119}
\begin{split}
&h_{\beta_{1},1}=\tilde{h}_{\beta_{1},1}-\frac{\hbar}{2}h_{\beta_{1},0}^{2}\rightarrow \tilde{h}_{\alpha_1}+\tilde{h}_{\alpha_2}+\frac{\hbar}{2}\{x_{\alpha_{2},0}^{+},x_{\alpha_{2},0}^{-}\} - \frac{\hbar}{2}(h_{\alpha_{1},0}+h_{\alpha_{2},0})^{2}=\\
&=h_{\alpha_{1},1}+h_{\alpha_{2},1}+ \frac{\hbar}{2}\{x_{\alpha_{2},0}^{+},x_{\alpha_{2},0}^{-}\} - \frac{\hbar}{2}\{h_{\alpha_{1},0},h_{\alpha_{2},0}\}.
\end{split}
\end{equation}
Thus by combining (\ref{equ:116}), (\ref{equ:117}), (\ref{equ:118}) and (\ref{equ:119}) we obtain that relation (\ref{rel:36}) holds. Thus we proved the Lemma (\ref{lemma:4.11}).

The inverse map $\psi^{-1}$ is map induced by simple reflection $s_{-\alpha_i}$. The proof of the relations in $\psi^{-1}$ case is similar to $\psi$ case, hence we obtained isomorphism map for $T_{s_{\alpha_i}}$, since every element of Weyl groupoid can be presented as product of elements induced by simple reflections, we can construct sequence of isomorphism maps for every element of groupoid. Thus we proved the statement.

\section{Conclusion}

Two presentations, a minimalistic presentation and a Drinfeld presentation of the Yangian of an affine Lie superalgebra $\hat{sl}(\Pi) = \hat{sl}(m|n)$ given by an arbitrary system of simple roots $\Pi$, are obtained in this work. It is proved that these two presentations are equivalent, that is, the Yangians given by these presentations are isomorphic. It is proved that the Yangians of affine Lie superalgebras $\hat{sl}(\Pi_1)$ and $\hat{sl}(\Pi_2)$  defined by arbitrary different systems of simple roots $\Pi_1$ and $\Pi_2$ are also isomorphic. We also give a construction of a Yangian Weyl groupoid whose elements define isomorphisms between super Yangians of affine Lie superalgebras defined by systems of simple roots of type $A$. The question of studying comultiplication structures on the affine Yangian is beyond the scope of this work. It is expected that it will also be investigated.\\

This work is partially supported by the Moscow Institute of Physics and Technology under the Priority 2030 Strategic Academic Leadership Program and by Russian Science Foundation grant 23-21-00282.

This work is performed at the Center of Pure Mathematics, MIPT, with
financial support of the project FSMG-2023-0013.



\begin{thebibliography}{00}
\bibitem{Dr} V.G. Drinfeld, Quantum groups, J. Math. Sci. 41 (1988) 898-915.
\bibitem{Levendorski}
S.I. Boyarchenko and S.Z. Levendorskii,  On affine Yangians.  Lett. Math. Phys., V. 32 (1993), no 4, 2691--274
\bibitem{Guaywork1}
N. Guay. Affine Yangians and deformed double current algebras in type A. Adv. Math., 211(2):436–484, 2007.
\bibitem{Guaywork}
H. Nakajima, N. Guay and C. Wendlandt. Coproduct for Yangians of affine Kac-Moody algebras on generators and
defining relations of Yangians. Adv. Math., 338:865–911, 2018.
\bibitem{Musson}
I. M. Musson. Lie superalgebras and enveloping algebras. Volume 131 of Graduate Studies in Mathematics., 2012.
\bibitem{Drinfeld}
V. G. Drinfeld. A new realization of Yangians and of quantum affine algebras. Dokl. Akad. Nauk SSSR, 296(1):,
pages 13–17, 1987
\bibitem{Lev} S.Z. Levendorskii, On generators and defining relations of Yangians. J. Geom. Phys. 12 (1993), no. 1, 1 -- 11.
\bibitem{SuSy}
C. Peng M. R. Gaberdiel, W. Li and H. Zhang. The supersymmetric affine Yangian,. 2018.
\bibitem{Ueda}
Mamoru Ueda. Construction of affine super Yangian. arXiv:1911.06666, 2021.
\bibitem{Kodera}
Ryosuke Kodera. Braid group action on affine Yangian. arXiv:1805.01621, 2019.
\bibitem{Mazurenko}
Vladimir A. Stukopin, Alexander Mazurenko. Classification of Hopf superalgebra structures on Drinfeld super
Yangians. arxiv:2210.08365, 2022.
\bibitem{Peng}
Y. N. Peng. On shifted super Yangians and a class of finite W-superalgebras. J. Algebra, 422, p. 520 -- 562, 2015.
\bibitem{PS}
Elena Poletaeva, Vera Serganova
 Representations of principal W-algebra for the superalgebra Q(n) and the super Yangian YQ(1). arXiv:1903.05272, 2019.
\bibitem{RS}
E. Ragoucy and P. Sorba. Yangian realisations from finite W-algebras. Comm. Math. Phys., 203(3), p. 551 -– 572, 1999.
\bibitem{St} V. Stukopin, Yangians of Lie superalgebras of type $A(m, n)$, (Russian) Funktsional. Anal. i Prilozhen. 28 (1994), no 3, 85 -- 88; translation in Funct. Anal. Appl. 28 (1994), no. 3, 217 -- 219.
\end{thebibliography}
\end{document}